# Fokas-type closed-form solution formulae for Sobolev-type equations with time-dependent coefficients

Andreas Chatziafratis

Division of Applied Analysis & PDE, Department of Mathematics, National & Kapodistrian University of Athens

December 17, 2025

**Abstract.** We analytically derive novel explicit integral representations for the solution of nonhomogeneous initial-boundary-value problems for a large category of evolution partial differential equations of Sobolev-Galpern type with generic temporally variable coefficients, satisfying suitable mild conditions, and with arbitrary data in classical function spaces. This work is based on the careful implementation of the pioneering Fokas unified transform methodology alongside its recently-proposed extension for solving a class of linear evolution equations with dispersion relation of specific polynomial type and time-dependent coefficients. We herein effectively extend those techniques to a special collection of evolution equations with time-dependent coefficients and mixed spatiotemporal derivatives, which induce rational dispersion relations. The new approach is exhibited in detail through illustrative generation of closed-form solutions for a multitude of such equations (such as Milne-Taylor-Barenblatt-Coleman-Ting-Chen-type, Benjamin-Bona-Mahony-type, as well as numerous higher-order variants) posed on the half-line. Challenging technical difficulties of complex-analytic and of algebraic flavour naturally emerge due the presence of mixed-derivative terms, and these are appropriately resolved in each case. The new formulas are of utility in subsequent investigation of qualitative properties and analysis of nonlinear counterparts too. Further extensions, generalizations, rigorous aspects and implementations to other types of problems as well will soon be reported in forthcoming publications.

## 1. Introduction

In what follows, addressed are initial-boundary-value problems (IBVP) posed on the half-line for a class of PDE with mixed derivatives (Sobolev-Galpern type [1]) and with general time-dependent coefficients and data, under mild appropriate assumptions.

We proceed in the spirit of the groundbreaking Fokas UTM methodology (see e.g. [2-9]). Inspiration is partly drawn from a beautiful observation reported in the recent work of Kalimeris and Özsarı who, in [10], indicated the way forward to producing, in a simple and elegant manner, UTM-type solution formulas for certain "polynomial-type" evolution PDE with time-variable coefficients. Motivation for studying PDE with such coefficients is also nicely provided in [10]. (It is worth mentioning that linear-PDE problems with spatially varying coefficients have been investigated via the Fokas method also in [11-14].)

The above approaches and ideas are here extended towards a new, significant direction. More precisely we will deal with equations of the form
$$\mathcal{P}u = f(x,t),$$
where $\mathcal{P}$ is one of the following differential operators:

1. $\partial_t - \alpha(t)\partial_{xxt} - \beta(t)\partial_{xx} + \gamma(t)$,
2. $\partial_t - \alpha\partial_{xxt} + \beta(t)\partial_x$,
3. $\partial_t + \alpha(t)\partial_{xxxxt} + \beta(t)\partial_{xxxx} + \gamma(t)$,
4. $\partial_t + \alpha\partial_{xxxxt} - \beta(t)\partial_{xx}$,
5. $\partial_t + \alpha\partial_{xxxxt} - \beta\partial_{xxt} + \gamma(t)\partial_{xxxx}$,
6. $\partial_t + \alpha\partial_{xxxxt} - \beta\partial_{xxt} - \gamma(t)\partial_x$,
7. $\partial_t - \alpha\partial_{xxxxxxt} + \beta\partial_{xxxxt} - \gamma\partial_{xxt} - \delta(t)\partial_{xx}$,
8. $\partial_t - \alpha(t)\partial_{xxxxxxt} - \beta(t)\partial_{xxxxxx} + \gamma(t)$,
9. $\partial_t + (-1)^\nu \alpha(t)\partial_x^{2\nu}\partial_t + (-1)^\nu \beta(t)\partial_x^{2\nu} + \gamma(t)$, $\nu \in \mathbb{N}$.

---

*E-mail address:* chatziafrati@math.uoa.gr

A general pattern is thence revealed from our demonstration which ultimately concerns a broad class of related problems. Such PDE and IBVP already have been associated with important applications in miscellaneous areas of applied sciences and engineering and further practical relevance is expected to be uncovered for fields such as solid mechanics, fluid dynamics, heat and mass transport, materials science, biology, mathematical finance, soil physics and nanotechnology.

Indicatively, let us observe that the 1[st] operator case alone (even with the coefficients replaced by constants and the reaction and forcing terms dropped) corresponds to the celebrated "pseudo-parabolic" model which is of paramount importance and has captured extensive attention in a large body of applied mathematics, continuum mechanics and mathematical physics literature. More specifically, such linear pseudo-parabolic equations have been derived or emerged, and studied from various aspects, in, e.g., [15-35]. Nevertheless, exact solution formulas for cases of time-dependent coefficients have apparently never been proposed systematically.

Further rigorous considerations, qualitative theory, asymptotics, controllability properties, extensions to non-standard problems, as well as applications to the analysis of nonlinear analogues will soon appear separately.

*Notation and assumptions* We will use the standard notation for the half-line Fourier transform. Thus, for $\lambda \in \mathbb{C}$ with $\operatorname{Im}\lambda \leq 0$, we define

$$\hat{u} = \hat{u}(\lambda,t) = \int_{y=0}^{\infty} u(y,t) e^{-i\lambda y} dy \text{ and } \hat{f}(\lambda,t) = \int_{y=0}^{\infty} f(y,t) e^{-i\lambda y} dy.$$

Throughtout this paper, we make the following assumptions on the data:

$$u_0(x) \in \mathcal{S}([0,\infty)), \ g_0(t), g_1(t), g_2(t), \ldots \in C^{\infty}([0,\infty)) \text{ and}$$

$$f = f(x,t) \in C^{\infty}(\overline{\mathbb{R}^+ \times \mathbb{R}^+}) \text{ such that } f(\cdot,t) \in \mathcal{S}([0,\infty)).$$

More precisely, the last assumption on the function $f(x,t)$ means that it is rapidly decreasing with respect to $x$, uniformy for $t$ in compact subsets of $[0,+\infty)$, i.e., for every $\ell, n \in \mathbb{N} \cup \{0\}$ and $t_0 > 0$,

$$\sup\left\{ x^{\ell} \left|\frac{\partial^n f(x,t)}{\partial x^n}\right| : x \geq 0, 0 \leq t \leq t_0 \right\} < +\infty.$$

Also, we use the following notation for the boundary values of the $x$–derivatives of the unknown function $u(x,t)$:

$$g_0(t) = u(0,t), \ g_1(t) = u_x(0,t), \ g_2(t) = u_{xx}(0,t), \ldots$$

Depending on the particular problem, the above functions $g_j(t)$ either are part of the data or are unknowns to be determined in the process of the contruction of the solution.

As customary, the derivation starts by assumimg that we have a solution $u(x,t)$ of the particular problem, which is appropriately behaved so that the steps and the computations of contruction are justified.

Regarding the coefficients $\alpha(t)$, $\beta(t)$, $\gamma(t)$, …, they are assumed to be $C^{\infty}$ for $t \in [0,+\infty)$, such that

$$\alpha(t) > 0, \text{ for every } t \in [0,+\infty) \text{ and } \sup\{\alpha(t): t \geq 0\} < +\infty.$$

Depending on the particular problem, further restrictions may have to be imposed. For example, in section 2, we assume $\beta(t) \not\equiv 0$. In section 6, we assume that $\alpha(t)$ and $\beta(t)$ are constants.

Furthermore, certain compatibility conditions are required. For example, in section 2, we assume $u_0(0) = g_0(0)$.

## 2. The operator $\partial_t - \alpha(t)\partial_{xxt} - \beta(t)\partial_{xx} + \gamma(t)$

**Problem 1** Solve

(2.1)
$$\begin{cases} \partial_t u = \alpha(t)\partial_{xxt}u + \beta(t)\partial_{xx}u - \gamma(t)u + f(x,t), \ (x,t) \in \mathbb{R}^+ \times \mathbb{R}^+, \\ u(x,0) = u_0(x), \ x \in \mathbb{R}^+, \\ u(0,t) = g_0(t), \ t \in \mathbb{R}^+, \end{cases}$$

for $u(x,t)$.

**Derivation of the solution** Working with $\lambda \in \mathbb{C}$, $\mathrm{Im}\,\lambda \leq 0$, and taking Fourier transforms, the diferential equation in (2.1) gives

(2.2) $\quad \dfrac{\partial}{\partial t}[\hat{u}(\lambda,t)] = \alpha(t)\dfrac{\partial}{\partial t}\left[\int_0^\infty \dfrac{\partial^2 u(y,t)}{\partial y^2}e^{-i\lambda y}dy\right] + \beta(t)\int_0^\infty \dfrac{\partial^2 u(y,t)}{\partial y^2}e^{-i\lambda y}dy - \gamma(t)\int_0^\infty u(y,t)e^{-i\lambda y}dy$,

since

$$\hat{u}_t = \hat{u}_t(\lambda,t) = \dfrac{\partial \hat{u}(\lambda,t)}{\partial t} = \int_{y=0}^\infty \dfrac{\partial u(y,t)}{\partial t}e^{-i\lambda y}dy.$$

Also,

$$\int_0^\infty u_x(x,t)e^{-i\lambda x}dx = -g_0(t) + i\lambda\hat{u}(\lambda,t) \text{ and } \int_0^\infty u_{xx}(x,t)e^{-i\lambda x}dx = -g_1(t) - i\lambda g_0(t) + (i\lambda)^2 \hat{u}(\lambda,t),$$

where we have set $g_1(t) = u_x(0,t)$, so, integrating by parts in (2.2), we obtain:

$$[1 + \alpha(t)\lambda^2]\dfrac{\partial}{\partial t}[\hat{u}(\lambda,t)] + [\beta(t)\lambda^2 + \gamma(t)]\hat{u}(\lambda,t) = -\alpha(t)[g_1'(t) + i\lambda g_0'(t)] - \beta(t)[g_1(t) + i\lambda g_0(t)] + \hat{f}(\lambda,t),$$

since $\dfrac{\partial}{\partial t}[u_x(0,t)] = g_1'(t)$.

Thus,

(2.3) $\quad \dfrac{\partial}{\partial t}[\hat{u}(\lambda,t)] + \dfrac{\beta(t)\lambda^2 + \gamma(t)}{1+\alpha(t)\lambda^2}\hat{u}(\lambda,t) = -\dfrac{1}{1+\alpha(t)\lambda^2}\{\alpha(t)[g_1'(t) + i\lambda g_0'(t)] + \beta(t)[g_1(t) + i\lambda g_0(t)]\} + \dfrac{\hat{f}(\lambda,t)}{1+\alpha(t)\lambda^2}.$

Setting

$$\omega(\lambda,t) = \dfrac{\beta(t)\lambda^2 + \gamma(t)}{1+\alpha(t)\lambda^2} \text{ and } \Omega(\lambda,t) = \int_{\tau=0}^t \omega(\lambda,\tau)d\tau = \int_{\tau=0}^t \dfrac{\beta(\tau)\lambda^2 + \gamma(\tau)}{1+\alpha(\tau)\lambda^2}d\tau,$$

we see – and this is crucial for what follows – that

(2.4) $\quad \omega(-\lambda,t) = \omega(\lambda,t), \ \Omega(-\lambda,t) = \Omega(\lambda,t),$

and

$$\lim_{\substack{\lambda \to \infty \\ \lambda \in \mathbb{C}-i\mathbb{R}}} \Omega(\lambda,t) = \int_{\tau=0}^t \dfrac{\beta(\tau)}{\alpha(\tau)}d\tau, \text{ uniformly for } t \text{ in compact subsets of } [0, +\infty).$$

Also, it follows from the above relation that, for every $T > 0$ and for every neighborhood A of the half-line $[c_0 i, +\infty i) := \{\eta i : \eta \geq c_0\}$, where $c_0 := \inf\{1/\sqrt{\alpha(t)} : t \geq 0\} > 0\} < +\infty$,

(2.5) $\quad \sup\{|e^{\pm\Omega(\lambda,t)}| : 0 \leq t \leq T, \ \lambda \in \mathbb{C} - \mathrm{A} \text{ and } \mathrm{Im}\,\lambda \geq 0\} < +\infty.$

Since

$$\frac{\partial}{\partial t}[e^{\Omega(\lambda,t)}\hat{u}(\lambda,t)] = e^{\Omega(\lambda,t)}\frac{\partial}{\partial t}[\hat{u}(\lambda,t)] + e^{\Omega(\lambda,t)}\frac{\beta(t)\lambda^2 + \gamma(t)}{1+\alpha(t)\lambda^2}\hat{u}(\lambda,t),$$

we may write (2.3) as follows:

$$\frac{\partial}{\partial t}[e^{\Omega(\lambda,t)}\hat{u}(\lambda,t)] = -\frac{e^{\Omega(\lambda,t)}}{1+\alpha(t)\lambda^2}\{\alpha(t)[g_1'(t)+i\lambda g_0'(t)] + \beta(t)[g_1(t)+i\lambda g_0(t)]\} + \frac{e^{\Omega(\lambda,t)}}{1+\alpha(t)\lambda^2}\hat{f}(\lambda,t).$$

Integrating the above equation, we obtain

$$e^{\Omega(\lambda,t)}\hat{u}(\lambda,t) = \hat{u}_0(\lambda) - \int_{\tau=0}^{t}\frac{e^{\Omega(\lambda,\tau)}}{1+\alpha(\tau)\lambda^2}\left\{\alpha(\tau)[g_1'(\tau)+i\lambda g_0'(\tau)] + \beta(\tau)[g_1(\tau)+i\lambda g_0(\tau)] - \hat{f}(\lambda,\tau)\right\}d\tau,$$

or equivalently,

$$(2.6) \quad \hat{u}(\lambda,t) = e^{-\Omega(\lambda,t)}\hat{u}_0(\lambda) - e^{-\Omega(\lambda,t)}\int_{\tau=0}^{t}\frac{e^{\Omega(\lambda,\tau)}}{1+\alpha(\tau)\lambda^2}[\alpha(\tau)g_1'(\tau) + \beta(\tau)g_1(\tau)]d\tau$$

$$-e^{-\Omega(\lambda,t)}\int_{\tau=0}^{t}\frac{e^{\Omega(\lambda,\tau)}}{1+\alpha(\tau)\lambda^2}\left\{i\lambda[\alpha(\tau)g_0'(\tau) + \beta(\tau)g_0(\tau)] - \hat{f}(\lambda,\tau)\right\}d\tau, \text{ for } \lambda \in \mathbb{C} \text{ with } \text{Im}\lambda \leq 0.$$

In view of (2.4), setting $-\lambda$ in place of $\lambda$, in (2.6), we find

$$(2.7) \quad \hat{u}(-\lambda,t) = e^{-\Omega(\lambda,t)}\hat{u}_0(-\lambda) - e^{-\Omega(\lambda,t)}\int_{\tau=0}^{t}\frac{e^{\Omega(\lambda,\tau)}}{1+\alpha(\tau)\lambda^2}[\alpha(\tau)g_1'(\tau) + \beta(\tau)g_1(\tau)]d\tau$$

$$-e^{-\Omega(\lambda,t)}\int_{\tau=0}^{t}\frac{e^{\Omega(\lambda,\tau)}}{1+\alpha(\tau)\lambda^2}\left\{-i\lambda[\alpha(\tau)g_0'(\tau) + \beta(\tau)g_0(\tau)] - \hat{f}(-\lambda,\tau)\right\}d\tau, \text{ for } \lambda \in \mathbb{C} \text{ with } \text{Im}\lambda \geq 0.$$

Multiplying (2.6) by $e^{i\lambda x}$ and integrating with respect to $\lambda$, we obtain

$$(2.8) \quad 2\pi u(x,t) = \int_{-\infty}^{\infty}e^{i\lambda x}\hat{u}(\lambda,t)d\lambda$$

$$= \int_{-\infty}^{\infty}e^{i\lambda x - \Omega(\lambda,t)}\hat{u}_0(\lambda)d\lambda - \int_{-\infty}^{\infty}\left\{e^{i\lambda x - \Omega(\lambda,t)}\int_{\tau=0}^{t}\frac{e^{\Omega(\lambda,\tau)}}{1+\alpha(\tau)\lambda^2}[\alpha(\tau)g_1'(\tau) + \beta(\tau)g_1(\tau)]d\tau\right\}d\lambda$$

$$-\int_{-\infty}^{\infty}\left\{i\lambda e^{i\lambda x - \Omega(\lambda,t)}\int_{\tau=0}^{t}\frac{e^{\Omega(\lambda,\tau)}}{1+\alpha(\tau)\lambda^2}[\alpha(\tau)g_0'(\tau) + \beta(\tau)g_0(\tau)]d\tau\right\}d\lambda + \int_{-\infty}^{\infty}\left\{e^{i\lambda x - \Omega(\lambda,t)}\int_{\tau=0}^{t}\frac{e^{\Omega(\lambda,\tau)}}{1+\alpha(\tau)\lambda^2}\hat{f}(\lambda,\tau)\,d\tau\right\}d\lambda.$$

Using (2.5) and deforming the contours of the $d\lambda$-integrals which involve the quantities

$$\alpha(\tau)g_1'(\tau) + \beta(\tau)g_1(\tau), \ \alpha(\tau)g_1'(\tau) + \beta(\tau)g_1(\tau),$$

we write (2.8) in the following way:

$$(2.8) \quad 2\pi u(x,t) = \int_{-\infty}^{\infty}e^{i\lambda x}\hat{u}(\lambda,t)d\lambda$$

$$= \int_{-\infty}^{\infty}e^{i\lambda x - \Omega(\lambda,t)}\hat{u}_0(\lambda)d\lambda - \int_{L}\left\{e^{i\lambda x - \Omega(\lambda,t)}\int_{\tau=0}^{t}\frac{e^{\Omega(\lambda,\tau)}}{1+\alpha(\tau)\lambda^2}[\alpha(\tau)g_1'(\tau) + \beta(\tau)g_1(\tau)]d\tau\right\}d\lambda$$

$$-\int_{L}\left\{i\lambda e^{i\lambda x - \Omega(\lambda,t)}\int_{\tau=0}^{t}\frac{e^{\Omega(\lambda,\tau)}}{1+\alpha(\tau)\lambda^2}[\alpha(\tau)g_0'(\tau) + \beta(\tau)g_0(\tau)]d\tau\right\}d\lambda + \int_{-\infty}^{\infty}\left\{e^{i\lambda x - \Omega(\lambda,t)}\int_{\tau=0}^{t}\frac{e^{\Omega(\lambda,\tau)}}{1+\alpha(\tau)\lambda^2}\hat{f}(\lambda,\tau)\,d\tau\right\}d\lambda,$$

where $\mathcal{L}$ is a contour as the one depicted in fig.1, below. More precisely it consists of two half-lines, parallel to the imaginary axis, and a half-circle-like curve, joining the half-lines, and oriented counterclockwise, as in fig. 1. "Minor" deformations of the half-lines will not change the values of the integrals in (2.8), which are taken over such a deformed $\mathcal{L}$.

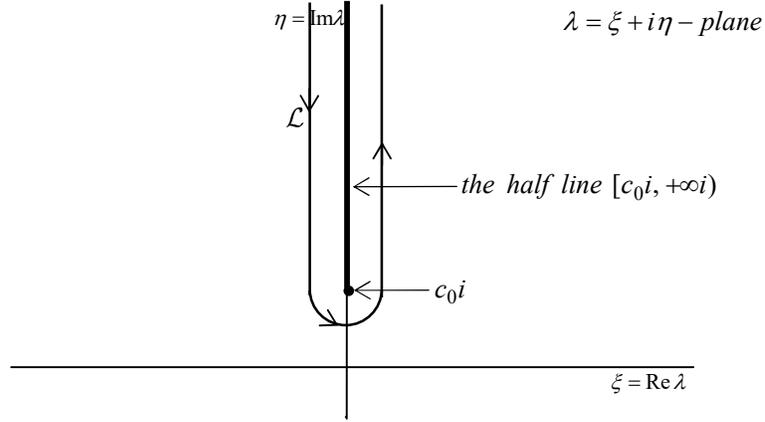

Fig. 1 A choice of the contour $\mathcal{L}$.

On the other hand, (2.7) gives

$$(2.9) \quad 0 = \int_{\mathcal{L}} e^{i\lambda x} \hat{u}(-\lambda, t) d\lambda$$

$$= \int_{\mathcal{L}} e^{i\lambda x - \Omega(\lambda,t)} \hat{u}_0(-\lambda) d\lambda - \int_{\mathcal{L}} \left\{ e^{i\lambda x - \Omega(\lambda,t)} \int_{\tau=0}^{t} \frac{e^{\Omega(\lambda,\tau)}}{1+\alpha(\tau)\lambda^2} [\alpha(\tau) g_1'(\tau) + \beta(\tau) g_1(\tau)] d\tau \right\} d\lambda$$

$$+ \int_{\mathcal{L}} \left\{ i\lambda e^{i\lambda x - \Omega(\lambda,t)} \int_{\tau=0}^{t} \frac{e^{\Omega(\lambda,\tau)}}{1+\alpha(\tau)\lambda^2} [\alpha(\tau) g_0'(\tau) + \beta(\tau) g_0(\tau)] d\tau \right\} d\lambda + \int_{\mathcal{L}} \left\{ e^{i\lambda x - \Omega(\lambda,t)} \int_{\tau=0}^{t} \frac{e^{\Omega(\lambda,\tau)}}{1+\alpha(\tau)\lambda^2} \hat{f}(-\lambda,\tau) d\tau \right\} d\lambda.$$

Subtracting (2.9) from (2.8), we obtain the solution to the Dirichlet problem 1: For $x > 0$ and $t > 0$,

$$(2.10) \quad 2\pi u(x,t) = \int_{-\infty}^{\infty} e^{i\lambda x - \Omega(\lambda,t)} \hat{u}_0(\lambda) d\lambda - \int_{\mathcal{L}} e^{i\lambda x - \Omega(\lambda,t)} \hat{u}_0(-\lambda) d\lambda$$

$$- 2\int_{\mathcal{L}} \left\{ i\lambda e^{i\lambda x - \Omega(\lambda,t)} \int_{\tau=0}^{t} \frac{e^{\Omega(\lambda,\tau)}}{1+\alpha(\tau)\lambda^2} [\alpha(\tau) g_0'(\tau) + \beta(\tau) g_0(\tau)] d\tau \right\} d\lambda$$

$$+ \int_{-\infty}^{\infty} \left\{ e^{i\lambda x - \Omega(\lambda,t)} \int_{\tau=0}^{t} \frac{e^{\Omega(\lambda,\tau)}}{1+\alpha(\tau)\lambda^2} \hat{f}(\lambda,\tau) d\tau \right\} d\lambda - \int_{\mathcal{L}} \left\{ e^{i\lambda x - \Omega(\lambda,t)} \int_{\tau=0}^{t} \frac{e^{\Omega(\lambda,\tau)}}{1+\alpha(\tau)\lambda^2} \hat{f}(-\lambda,\tau) d\tau \right\} d\lambda.$$

*Remark* In the case the function $\alpha(t)$ satisfies, in addition, the condition

$$C_0 := \sup\{1/\sqrt{\alpha(t)} : t \geq 0\} < +\infty,$$

the contour $\mathcal{L}$, in formula (2.10), can be replaced by any finite closed contour $\mathcal{C}$, surrounding the line segment $[c_0 i, C_0 i] = \{\eta i : c_0 \leq \eta \leq C_0\}$ and oriented counterclockwise, as depicted in fig.2, below.

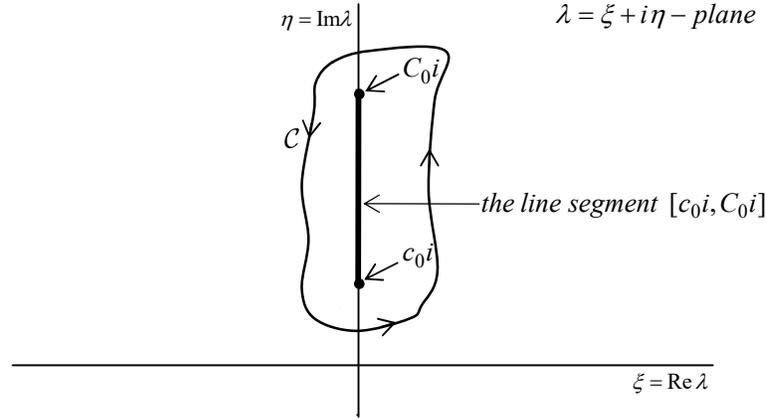

Fig. 2 A choice of the contour $\mathcal{C}$.

***The Neumann problem*** *Solve*

(2.11)
$$\begin{cases} \partial_t u = \alpha(t)\partial_{xxt} u + \beta(t)\partial_{xx} u - \gamma(t)u + f(x,t),\ (x,t) \in \mathbb{R}^+ \times \mathbb{R}^+, \\ u(x,0) = u_0(x),\ x \in \mathbb{R}^+, \\ u_x(0,t) = g_1(t),\ t \in \mathbb{R}^+, \end{cases}$$

*for* $u(x,t)$.

Working as in the case of the Dirichlet problem, we find that the solution to the above Neumann problem is given by the following formula:

(2.12) $2\pi u(x,t) = \int_{-\infty}^{\infty} e^{i\lambda x - \Omega(\lambda,t)} \hat{u}_0(\lambda) d\lambda - \int_{\mathcal{L}} e^{i\lambda x - \Omega(\lambda,t)} \hat{u}_0(-\lambda) d\lambda$

$\quad - 2\int_{\mathcal{L}} \left\{ e^{i\lambda x - \Omega(\lambda,t)} \int_{\tau=0}^{t} \frac{e^{\Omega(\lambda,\tau)}}{1+\alpha(\tau)\lambda^2}[\alpha(\tau)g_1'(\tau) + \beta(\tau)g_1(\tau)] d\tau \right\} d\lambda$

$\quad + \int_{-\infty}^{\infty} \left\{ e^{i\lambda x - \Omega(\lambda,t)} \int_{\tau=0}^{t} \frac{e^{\Omega(\lambda,\tau)}}{1+\alpha(\tau)\lambda^2} \hat{f}(\lambda,\tau) d\tau \right\} d\lambda - \int_{\mathcal{L}} \left\{ e^{i\lambda x - \Omega(\lambda,t)} \int_{\tau=0}^{t} \frac{e^{\Omega(\lambda,\tau)}}{1+\alpha(\tau)\lambda^2} \hat{f}(-\lambda,\tau) d\tau \right\} d\lambda.$

## 3. The operator $\partial_t - \alpha \partial_{xxt} + \beta(t)\partial_x$

***Problem 2*** *Solve*

(3.1)
$$\begin{cases} \partial_t u = \alpha \partial_{xxt} u - \beta(t)\partial_x u + f(x,t),\ (x,t) \in Q, \\ u(x,0) = u_0(x),\ x \in \mathbb{R}^+, \\ u(0,t) = g_0(t),\ t \in \mathbb{R}^+, \end{cases}$$

*for* $u(x,t)$.

***Derivation of the solution*** Taking Fourier transforms, the diferential equation in (3.1) gives

(3.2) $\quad \dfrac{\partial}{\partial t}[\hat{u}(\lambda,t)] = \alpha \dfrac{\partial}{\partial t}\left[\int_0^{\infty} \dfrac{\partial^2 u(y,t)}{\partial y^2} e^{-i\lambda y} dy\right] - \beta(t)\int_0^{\infty} \dfrac{\partial u(y,t)}{\partial y} e^{-i\lambda y} dy.$

Integrating by parts in (3.2), we obtain:

$$(1+\alpha\lambda^2)\dfrac{\partial}{\partial t}[\hat{u}(\lambda,t)] + i\beta(t)\lambda \hat{u}(\lambda,t) = -\alpha[g_1'(t) + i\lambda g_0'(t)] + \beta(t)g_0(t) + \hat{f}(\lambda,t).$$

It follows that

$$\frac{\partial}{\partial t}[\hat{u}(\lambda,t)] + \frac{i\beta(t)\lambda}{1+\alpha\lambda^2}\hat{u}(\lambda,t) = \frac{1}{1+\alpha\lambda^2}[\beta(t)g_0(t) - i\alpha\lambda g_0'(t)] - \frac{\alpha}{1+\alpha\lambda^2}g_1'(t) + \frac{\hat{f}(\lambda,t)}{1+\alpha\lambda^2}.$$

Thus, setting

$$\omega(\lambda,t) = \frac{i\beta(t)\lambda}{1+\alpha\lambda^2} \quad \text{and} \quad \Omega(\lambda,t) = \int_{\tau=0}^{t} \frac{i\beta(\tau)\lambda}{1+\alpha\lambda^2}d\tau,$$

we have

$$\frac{\partial}{\partial t}[e^{\Omega(\lambda,t)}\hat{u}(\lambda,t)] = \frac{e^{\Omega(\lambda,t)}}{1+\alpha\lambda^2}[\beta(t)g_0(t) - i\alpha\lambda g_0'(t)] - \frac{\alpha e^{\Omega(\lambda,t)}}{1+\alpha\lambda^2}g_1'(t) + \frac{e^{\Omega(\lambda,t)}\hat{f}(\lambda,t)}{1+\alpha\lambda^2}.$$

Therefore,

$$e^{\Omega(\lambda,t)}\hat{u}(\lambda,t) = \hat{u}_0(\lambda) + \int_{\tau=0}^{t}\left\{\frac{e^{\Omega(\lambda,t)}}{1+\alpha\lambda^2}[\beta(t)g_0(t) - i\alpha\lambda g_0'(t)] - \frac{\alpha e^{\Omega(\lambda,t)}}{1+\alpha\lambda^2}g_1'(t) + \frac{e^{\Omega(\lambda,t)}\hat{f}(\lambda,t)}{1+\alpha\lambda^2}\right\}d\tau,$$

and, equivalently,

$$(3.3) \quad \hat{u}(\lambda,t) = e^{-\Omega(\lambda,t)}\hat{u}_0(\lambda) + \frac{e^{-\Omega(\lambda,t)}}{1+\alpha\lambda^2}\int_{\tau=0}^{t}e^{\Omega(\lambda,\tau)}[\beta(\tau)g_0(\tau) - i\alpha\lambda g_0'(\tau)]d\tau$$

$$- \frac{\alpha e^{-\Omega(\lambda,t)}}{1+\alpha\lambda^2}\int_{\tau=0}^{t}e^{\Omega(\lambda,\tau)}g_1'(\tau)d\tau + \frac{e^{-\Omega(\lambda,t)}}{1+\alpha\lambda^2}\int_{\tau=0}^{t}e^{\Omega(\lambda,\tau)}\hat{f}(\lambda,\tau)d\tau.$$

Multiplying (3.3) by $e^{i\lambda x}$ and integrating we obtain

$$(3.4) \quad 2\pi u(x,t) = \int_{-\infty}^{\infty}e^{i\lambda x}\hat{u}(\lambda,t))d\lambda$$

$$= \int_{-\infty}^{\infty}e^{i\lambda x - \Omega(\lambda,t)}\hat{u}_0(\lambda)d\lambda + \int_{-\infty}^{\infty}\left\{\frac{e^{i\lambda x - \Omega(\lambda,t)}}{1+\alpha\lambda^2}\int_{\tau=0}^{t}e^{\Omega(\lambda,\tau)}[\beta(\tau)g_0(\tau) - i\alpha\lambda g_0'(\tau)]d\tau\right\}d\lambda$$

$$- \int_{-\infty}^{\infty}\left\{\frac{\alpha e^{i\lambda x - \Omega(\lambda,t)}}{1+\alpha\lambda^2}\int_{\tau=0}^{t}e^{\Omega(\lambda,\tau)}g_1'(\tau)d\tau\right\}d\lambda + \int_{-\infty}^{\infty}\left\{\frac{e^{i\lambda x - \Omega(\lambda,t)}}{1+\alpha\lambda^2}\int_{\tau=0}^{t}e^{\Omega(\lambda,\tau)}\hat{f}(\lambda,\tau)d\tau\right\}d\lambda.$$

Since

$$\lim_{\lambda\to\infty}\Omega(\lambda,t) = 0, \text{ uniformly for } t \text{ in compact subsets of } [0,+\infty),$$

we can deform the contours of the integrals in (3.4), which contain the terms

$$\beta(\tau)g_0(\tau) - i\alpha\lambda g_0'(\tau) \text{ and } g_1'(\tau),$$

obtaining the formula

$$(3.5) \quad 2\pi u(x,t) = \int_{-\infty}^{\infty}e^{i\lambda x}\hat{u}(\lambda,t))d\lambda$$

$$= \int_{-\infty}^{\infty}e^{i\lambda x - \Omega(\lambda,t)}\hat{u}_0(\lambda)d\lambda + \int_{C}\left\{\frac{e^{i\lambda x - \Omega(\lambda,t)}}{1+\alpha\lambda^2}\int_{\tau=0}^{t}e^{\Omega(\lambda,\tau)}[\beta(\tau)g_0(\tau) - i\alpha\lambda g_0'(\tau)]d\tau\right\}d\lambda$$

$$- \int_{C}\left\{\frac{\alpha e^{i\lambda x - \Omega(\lambda,t)}}{1+\alpha\lambda^2}\int_{\tau=0}^{t}e^{\Omega(\lambda,\tau)}g_1'(\tau)d\tau\right\}d\lambda + \int_{-\infty}^{\infty}\left\{\frac{e^{i\lambda x - \Omega(\lambda,t)}}{1+\alpha\lambda^2}\int_{\tau=0}^{t}e^{\Omega(\lambda,\tau)}\hat{f}(\lambda,\tau)d\tau\right\}d\lambda,$$

where the contour $C$ is any simple closed curve, sufficiently small, surrounding the point $\lambda = i/\sqrt{\alpha}$, lying in the upper half-plane $\{\lambda \in C : \text{Im}\,\lambda > 0\}$ and oriented counterclockwise (see fig. 3).

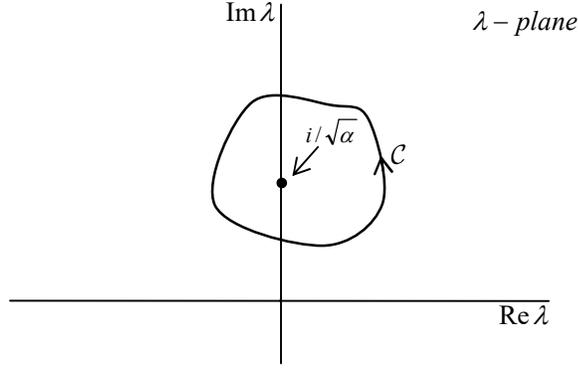

Fig. 3 A typical contour $\mathcal{C}$.

Next observing that

$$\omega(1/(\alpha\lambda),t) = \omega(\lambda,t) \text{ and } \Omega(1/(\alpha\lambda),t) = \Omega(\lambda,t),$$

and setting $1/(\alpha\lambda)$, in place of $\lambda$, in (3.3), we obtain

$$(3.6) \quad \hat{u}(1/(\alpha\lambda),t) = e^{-\Omega(\lambda,t)}\hat{u}_0(1/(\alpha\lambda)) + \frac{\alpha\lambda}{1+\alpha\lambda^2}e^{-\Omega(\lambda,t)}\int_{\tau=0}^{t}e^{\Omega(\lambda,\tau)}[\lambda\beta(\tau)g_0(\tau) - ig_0'(\tau)]d\tau$$

$$-\frac{\alpha\lambda^2}{1+\alpha\lambda^2}e^{-\Omega(\lambda,t)}\int_{\tau=0}^{t}e^{\Omega(\lambda,\tau)}g_1'(\tau)d\tau + \frac{e^{-\Omega(\lambda,t)}}{1+\alpha\lambda^2}\int_{\tau=0}^{t}e^{\Omega(\lambda,\tau)}\hat{f}(\lambda,\tau)d\tau.$$

Multiplying (3.6) by $e^{i\lambda x}/\lambda^2$ and integrating for $\lambda \in \mathcal{C}$, we see that

$$(3.7) \quad 0 = \int_{\mathcal{C}}\frac{1}{\lambda^2}e^{i\lambda x-\Omega(\lambda,t)}\hat{u}_0(1/(\alpha\lambda))d\lambda + \int_{\mathcal{C}}\left\{\frac{\alpha e^{i\lambda x-\Omega(\lambda,t)}}{\lambda(1+\alpha\lambda^2)}\int_{\tau=0}^{t}e^{\Omega(\lambda,\tau)}[\lambda\beta(\tau)g_0(\tau) - ig_0'(\tau)]d\tau\right\}d\lambda$$

$$-\int_{\mathcal{C}}\left\{\frac{\alpha e^{i\lambda x-\Omega(\lambda,t)}}{1+\alpha\lambda^2}\int_{\tau=0}^{t}e^{\Omega(\lambda,\tau)}g_1'(\tau)d\tau\right\}d\lambda + \int_{\mathcal{C}}\left\{\frac{e^{i\lambda x-\Omega(\lambda,t)}}{\lambda^2(1+\alpha\lambda^2)}\int_{\tau=0}^{t}e^{\Omega(\lambda,\tau)}\hat{f}(\lambda,\tau)d\tau\right\}d\lambda,$$

where we used the fact that

$$\int_{\mathcal{C}}\frac{e^{i\lambda x}}{\lambda^2}\hat{u}(1/(\alpha\lambda),t)d\lambda = 0.$$

Subtracting (3.7) from (3.5), we obtain the formula for the solution to the Dirichlet problem (3.1):

$$2\pi u(x,t) = \int_{-\infty}^{\infty}e^{i\lambda x-\Omega(\lambda,t)}\hat{u}_0(\lambda)d\lambda - \int_{\mathcal{C}}\frac{1}{\lambda^2}e^{i\lambda x-\Omega(\lambda,t)}\hat{u}_0(1/(\alpha\lambda))d\lambda$$

$$+\int_{\mathcal{C}}\left\{\frac{e^{i\lambda x-\Omega(\lambda,t)}}{1+\alpha\lambda^2}\int_{\tau=0}^{t}e^{\Omega(\lambda,\tau)}\left[(1-\alpha)\beta(\tau)g_0(\tau) + \frac{i\alpha(1-\lambda^2)}{\lambda}g_0'(\tau)\right]d\tau\right\}d\lambda$$

$$\int_{-\infty}^{\infty}\left\{\frac{e^{i\lambda x-\Omega(\lambda,t)}}{1+\alpha\lambda^2}\int_{\tau=0}^{t}e^{\Omega(\lambda,\tau)}\hat{f}(\lambda,\tau)d\tau\right\}d\lambda - \int_{\mathcal{C}}\left\{\frac{e^{i\lambda x-\Omega(\lambda,t)}}{\lambda^2(1+\alpha\lambda^2)}\int_{\tau=0}^{t}e^{\Omega(\lambda,\tau)}\hat{f}(\lambda,\tau)d\tau\right\}d\lambda.$$

## 4. The operator $\partial_t + \alpha(t)\partial_{xxxxt} + \beta(t)\partial_{xxxx} + \gamma(t)$

***Problem 3*** Solve

(4.1)
$$\begin{cases} \partial_t u = -\alpha(t)\partial_{xxxxt} u - \beta(t)\partial_{xxxx} u - \gamma(t)u + f(x,t), \ (x,t) \in \mathbb{R}^+ \times \mathbb{R}^+, \\ u(x,0) = u_0(x), \ x \in \mathbb{R}^+, \\ u(0,t) = g_0(t), \ t \in \mathbb{R}^+, \\ u_x(0,t) = g_1(t), \ t \in \mathbb{R}^+, \end{cases}$$

for $u(x,t)$.

***Derivation of the solution*** With notation as in the previous sections, taking Fourier transforms and integrating by parts, we see that the differential equation in (4.1) gives, for $\lambda \in \mathbb{C}$ with $\mathrm{Im}\,\lambda \leq 0$,

$$\frac{\partial}{\partial t}[\hat{u}(\lambda,t)] + \alpha(t)\left\{\lambda^4 \frac{\partial}{\partial t}[\hat{u}(\lambda,t)] - g'_3(t) - i\lambda g'_2(t) - (i\lambda)^2 g'_1(t) - (i\lambda)^3 g'_0(t)\right\}$$
$$+ \beta(t)[\lambda^4 \hat{u}(\lambda,t) - g_3(t) - i\lambda g_2(t) - (i\lambda)^2 g_1(t) - (i\lambda)^3 g_0(t)] + \gamma(t)\hat{u}(\lambda,t) = \hat{f}(\lambda,t),$$

$\Rightarrow$

$$[1 + \alpha(t)\lambda^4]\frac{\partial}{\partial t}[\hat{u}(\lambda,t)] + [\beta(t)\lambda^4 + \gamma(t)]\hat{u}(\lambda,t) = \alpha(t)[g'_3(t) + i\lambda g'_2(t) - \lambda^2 g'_1(t) - i\lambda^3 g'_0(t)]$$
$$+ \beta(t)[g_3(t) + i\lambda g_2(t) - \lambda^2 g_1(t) - i\lambda^3 g_0(t)] + \hat{f}(\lambda,t)$$

$\Rightarrow$

(4.2) $\quad \dfrac{\partial}{\partial t}[\hat{u}(\lambda,t)] + \dfrac{\beta(t)\lambda^4 + \gamma(t)}{1 + \alpha(t)\lambda^4}\hat{u}(\lambda,t) = \dfrac{1}{1+\alpha(t)\lambda^4}[G_3(t) + i\lambda G_2(t) - \lambda^2 G_1(t) - i\lambda^3 G_0(t)] + \dfrac{1}{1+\alpha(t)\lambda^4}\hat{f}(\lambda,t),$

where we have set

$$G_3(t) = \alpha(t)g'_3(t) + \beta(t)g_3(t), \ G_2(t) = \alpha(t)g'_2(t) + \beta(t)g_2(t),$$
$$G_1(t) = \alpha(t)g'_1(t) + \beta(t)g_1(t), \ G_0(t) = \alpha(t)g'_0(t) + \beta(t)g_0(t).$$

Setting,

$$\omega(\lambda,t) = \frac{\beta(t)\lambda^4 + \gamma(t)}{1+\alpha(t)\lambda^4} \ \text{ and } \ \Omega(\lambda,t) = \int_{\tau=0}^{t}\omega(\lambda,\tau)d\tau = \int_{\tau=0}^{t}\frac{\beta(\tau)\lambda^4 + \gamma(\tau)}{1+\alpha(\tau)\lambda^4}d\tau,$$

we write (4.2) in the form:

$$\frac{\partial}{\partial t}[e^{\Omega(\lambda,t)}\hat{u}(\lambda,t)] = \frac{e^{\Omega(\lambda,t)}}{1+\alpha(t)\lambda^4}[G_3(t) + i\lambda G_2(t) - \lambda^2 G_1(t) - i\lambda^3 G_0(t)] + \frac{e^{\Omega(\lambda,t)}}{1+\alpha(t)\lambda^4}\hat{f}(\lambda,t).$$

Integrating the above equation, we obtain that, for $\lambda \in \mathbb{C}$ with $\mathrm{Im}\,\lambda \leq 0$,

(4.3) $\ e^{\Omega(\lambda,t)}\hat{u}(\lambda,t) = \hat{u}_0(\lambda) + \displaystyle\int_{\tau=0}^{t}\frac{e^{\Omega(\lambda,\tau)}}{1+\alpha(\tau)\lambda^4}[G_3(\tau) + i\lambda G_2(\tau) - \lambda^2 G_1(\tau) - i\lambda^3 G_0(\tau)]d\tau + \int_{\tau=0}^{t}\frac{e^{\Omega(\lambda,\tau)}}{1+\alpha(\tau)\lambda^4}\hat{f}(\lambda,\tau)d\tau.$

Setting,

$$\widetilde{G}_j(\Omega,\lambda,t) = \widetilde{G}_j(\Omega(\lambda,t),\lambda,t) = \int_{\tau=0}^{t}\frac{e^{\Omega(\lambda,\tau)}}{1+\alpha(\tau)\lambda^4}G_j(\tau)d\tau, \ j = 1, 2, 3, 4,$$

and

$$\widetilde{\hat{f}}(\Omega,\lambda,t) = \widetilde{\hat{f}}(\Omega(\lambda,t),\lambda,t) = \int_{\tau=0}^{t}\frac{e^{\Omega(\lambda,\tau)}}{1+\alpha(\tau)\lambda^4}\hat{f}(\lambda,\tau)d\tau,$$

we write (4,3) as follows:

(4.4) $\hat{u}(\lambda,t) = e^{-\Omega(\lambda,t)}\hat{u}_0(\lambda) + e^{-\Omega(\lambda,t)}\widetilde{G}_3(\Omega(\lambda,t),\lambda,t) + i\lambda e^{-\Omega(\lambda,t)}\widetilde{G}_2(\Omega(\lambda,t),\lambda,t)$

$- \lambda^2 e^{-\Omega(\lambda,t)}\widetilde{G}_1(\Omega(\lambda,t),\lambda,t) - i\lambda^3 e^{-\Omega(\lambda,t)}\widetilde{G}_0(\Omega(\lambda,t),\lambda,t) + e^{-\Omega(\lambda,t)}\widetilde{\tilde{f}}(\Omega(\lambda,t),\lambda,t)$.

Next, we seek for the functions $\sigma(\lambda)$ for which

$$\omega(\sigma(\lambda),t) = \omega(\lambda,t) \implies \frac{\beta(t)[\sigma(\lambda)]^4 + \gamma(t)}{1+\alpha(t)[\sigma(\lambda)]^4} = \frac{\beta(t)\lambda^4 + \gamma(t)}{1+\alpha(t)\lambda^4}.$$

Solving the above algebraic equation we find

(4.5) $\qquad\qquad\qquad \sigma(\lambda) = \lambda$ or $\sigma(\lambda) = -\lambda$ or $\sigma(\lambda) = i\lambda$ or $\sigma(\lambda) = -i\lambda$.

(For the method to work, it crucial that all the functions $\sigma$, which satisfy the equation $\omega(\sigma,t) = \omega(\lambda,t)$, do not depend on $t$.)

By the choice of the functions (4.5), we have

$$\Omega(\lambda,t) = \Omega(-\lambda,t) = \Omega(i\lambda,t) = \Omega(-i\lambda,t),$$

and

$$\frac{e^{\Omega(\lambda,t)}}{1+\alpha(t)\lambda^4} = \frac{e^{\Omega(-\lambda,t)}}{1+\alpha(t)(-\lambda)^4} = \frac{e^{\Omega(i\lambda,t)}}{1+\alpha(t)(i\lambda)^4} = \frac{e^{\Omega(-i\lambda,t)}}{1+\alpha(t)(-i\lambda)^4}.$$

Therefore, setting $-\lambda$, $i\lambda$ and $-i\lambda$, in (4.4), in place of $\lambda$, we obtain

(4.5) $\hat{u}(-\lambda,t) = e^{-\Omega(\lambda,t)}\hat{u}_0(-\lambda) + e^{-\Omega(\lambda,t)}\widetilde{G}_3(\Omega(\lambda,t),\lambda,t) - i\lambda e^{-\Omega(\lambda,t)}\widetilde{G}_2(\Omega(\lambda,t),\lambda,t)$

$-\lambda^2 e^{-\Omega(\lambda,t)}\widetilde{G}_1(\Omega(\lambda,t),\lambda,t) + i\lambda^3 e^{-\Omega(\lambda,t)}\widetilde{G}_0(\Omega(\lambda,t),\lambda,t) + e^{-\Omega(\lambda,t)}\widetilde{\tilde{f}}(\Omega(\lambda,t),-\lambda,t)$, with $\mathrm{Im}\,\lambda \geq 0$,

(4.6) $\hat{u}(i\lambda,t) = e^{-\Omega(\lambda,t)}\hat{u}_0(i\lambda) + e^{-\Omega(\lambda,t)}\widetilde{G}_3(\Omega(\lambda,t),\lambda,t) - \lambda^2 e^{-\Omega(\lambda,t)}\widetilde{G}_2(\Omega(\lambda,t),\lambda,t)$

$+\lambda^2 e^{-\Omega(\lambda,t)}\widetilde{G}_1(\Omega(\lambda,t),\lambda,t) - \lambda^3 e^{-\Omega(\lambda,t)}\widetilde{G}_0(\Omega(\lambda,t),\lambda,t) + e^{-\Omega(\lambda,t)}\widetilde{\tilde{f}}(\Omega(\lambda,t),i\lambda,t)$, with $\mathrm{Im}(i\lambda) \leq 0$,

and

(4.7) $\hat{u}(-i\lambda,t) = e^{-\Omega(\lambda,t)}\hat{u}_0(-i\lambda) + e^{-\Omega(\lambda,t)}\widetilde{G}_3(\Omega(\lambda,t),\lambda,t) + \lambda^2 e^{-\Omega(\lambda,t)}\widetilde{G}_2(\Omega(\lambda,t),\lambda,t)$

$+\lambda^2 e^{-\Omega(\lambda,t)}\widetilde{G}_1(\Omega(\lambda,t),\lambda,t) - i\lambda^3 e^{-\Omega(\lambda,t)}\widetilde{G}_0(\Omega(\lambda,t),\lambda,t) + e^{-\Omega(\lambda,t)}\widetilde{\tilde{f}}(\Omega(\lambda,t),\lambda,t)$, with $\mathrm{Im}(-i\lambda) \leq 0$.

Also, multiplying (4.4) by $e^{i\lambda x}$ and integrating we obtain, for $x > 0$ and $t > 0$,

(4.8) $2\pi u(x,t) = \int_{-\infty}^{\infty} e^{i\lambda x}\hat{u}(\lambda,t)d\lambda$

$= \int_{-\infty}^{\infty} e^{i\lambda x-\Omega(\lambda,t)}\hat{u}_0(\lambda)d\lambda + \int_{-\infty}^{\infty} e^{i\lambda x-\Omega(\lambda,t)}\widetilde{G}_3(\Omega(\lambda,t),\lambda,t)d\lambda + \int_{-\infty}^{\infty} i\lambda e^{i\lambda x-\Omega(\lambda,t)}\widetilde{G}_2(\Omega(\lambda,t),\lambda,t)d\lambda$

$- \int_{-\infty}^{\infty} \lambda^2 e^{i\lambda x-\Omega(\lambda,t)}\widetilde{G}_1(\Omega(\lambda,t),\lambda,t)d\lambda - \int_{-\infty}^{\infty} i\lambda^3 e^{i\lambda x-\Omega(\lambda,t)}\widetilde{G}_0(\Omega(\lambda,t),\lambda,t)d\lambda + \int_{-\infty}^{\infty} e^{i\lambda x-\Omega(\lambda,t)}\widetilde{\tilde{f}}(\Omega(\lambda,t),\lambda,t)d\lambda$.

Next, observing that

$$1 + \alpha(t)\lambda^4 = 0 \implies \lambda = \pm e^{i\pi/4}/\sqrt[4]{\alpha(t)} \text{ or } \lambda = \pm e^{i3\pi/4}/\sqrt[4]{\alpha(t)},$$

we consider the functions

$$\lambda_1(t) := e^{i\pi/4}/\sqrt[4]{\alpha(t)} \text{ and } \lambda_2 := e^{i3\pi/4}/\sqrt[4]{\alpha(t)}.$$

Let us observe that

$$\lambda_1(t) \in \ell_1 := [c_0 e^{i\pi/4}, \infty e^{i\pi/4}) = \{c e^{i\pi/4} : c \geq c_0\} \text{ and } \lambda_2(t) \in \ell_2 := [c_0 e^{i3\pi/4}, \infty e^{i3\pi/4}) = \{c e^{i3\pi/4} : c \geq c_0\}.$$

(We assume that $c_0 := \inf\{1/\sqrt[4]{\alpha(t)} : t \geq 0\} > 0$.)

Fixing $x > 0$ and $t > 0$, and deforming the contours of the integrals which involve the quantities $\widetilde{G}_j(\Omega(\lambda,t),\lambda,t)$, we write (4.8) as follows:

$$(4.9) \quad 2\pi u(x,t) = \int_{-\infty}^{\infty} e^{i\lambda x - \Omega(\lambda,t)} \hat{u}_0(\lambda) d\lambda$$

$$+ \left(\int_{\mathcal{L}_1} + \int_{\mathcal{L}_2}\right)[e^{i\lambda x - \Omega(\lambda,t)} \widetilde{G}_3(\Omega(\lambda,t),\lambda,t) d\lambda] + \left(\int_{\mathcal{L}_1} + \int_{\mathcal{L}_2}\right)[i\lambda e^{i\lambda x - \Omega(\lambda,t)} \widetilde{G}_2(\Omega(\lambda,t),\lambda,t) d\lambda]$$

$$- \left(\int_{\mathcal{L}_1} + \int_{\mathcal{L}_2}\right)[\lambda^2 e^{i\lambda x - \Omega(\lambda,t)} \widetilde{G}_1(\Omega(\lambda,t),\lambda,t) d\lambda] - \left(\int_{\mathcal{L}_1} + \int_{\mathcal{L}_2}\right)[i\lambda^3 e^{i\lambda x - \Omega(\lambda,t)} \widetilde{G}_0(\Omega(\lambda,t),\lambda,t) d\lambda]$$

$$+ \int_{-\infty}^{\infty} e^{i\lambda x - \Omega(\lambda,t)} \widetilde{\widetilde{f}}(\Omega(\lambda,t),\lambda,t) d\lambda,$$

where the contours $\mathcal{L}_1$ and $\mathcal{L}_2$ are depicted in fig. 4, below. More precisely, $\mathcal{L}_1$ and $\mathcal{L}_2$ are the boundaries of sufficiently small, smoothly bounded, neighbourhoods of the half-lines $[c_0 e^{i\pi/4}, +\infty e^{i\pi/4})$ and $[c_0 e^{i3\pi/4}, +\infty e^{i3\pi/4})$, countclockwise oriented. (We assume that the width of the opening of the contours – at infinity – is bounded away from zero.)

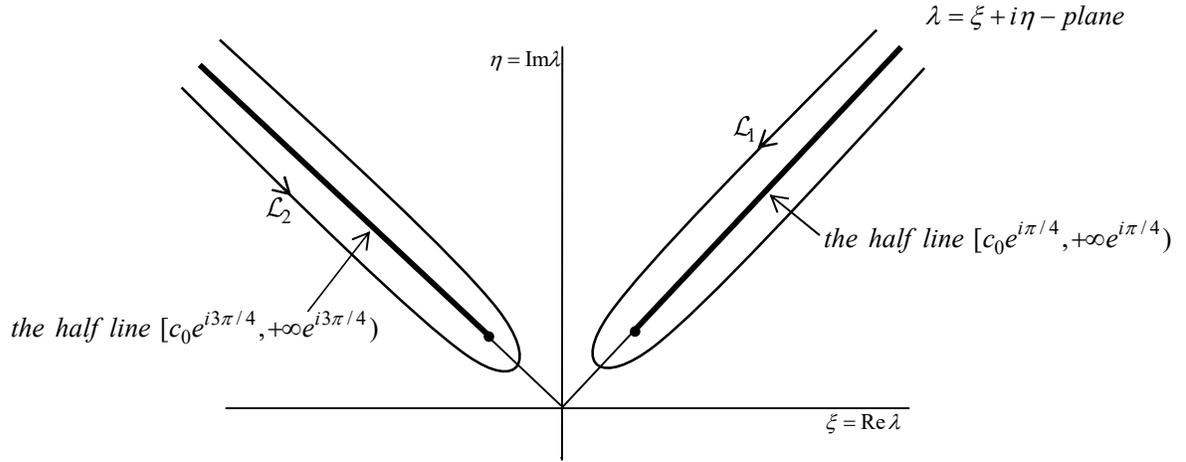

Fig. 4 Choices of the contours $\mathcal{L}_1$ and $\mathcal{L}_2$.

Using the abbreviations

$$\mathfrak{G}_j = \mathfrak{G}_j(\lambda,x,t) := e^{i\lambda x - \Omega(\lambda,t)} \widetilde{G}_j(\Omega(\lambda,t),\lambda,t), \quad j = 1, 2, 3, 4,$$

we write (4.9) as

$$(4.10) \quad 2\pi u(x,t) = \int_{-\infty}^{\infty} e^{i\lambda x - \Omega(\lambda,t)} \hat{u}_0(\lambda) d\lambda$$

$$+ \int_{\mathcal{L}_1} \mathfrak{G}_3(\lambda,x,t)d\lambda + \int_{\mathcal{L}_2} \mathfrak{G}_3(\lambda,x,t)d\lambda + \int_{\mathcal{L}_1} i\lambda \mathfrak{G}_2(\lambda,x,t)d\lambda + \int_{\mathcal{L}_2} i\lambda \mathfrak{G}_2(\lambda,x,t)d\lambda$$

$$+ \int_{\mathcal{L}_1} \lambda^2 \mathfrak{G}_1(\lambda,x,t)d\lambda + \int_{\mathcal{L}_2} \lambda^2 \mathfrak{G}_1(\lambda,x,t)d\lambda + \int_{\mathcal{L}_1} i\lambda^3 \mathfrak{G}_0(\lambda,x,t)d\lambda + \int_{\mathcal{L}_2} i\lambda^3 \mathfrak{G}_0(\lambda,x,t)d\lambda d\lambda$$

$$+ \int_{-\infty}^{\infty} e^{i\lambda x - \Omega(\lambda,t)} \widetilde{\tilde{f}}(\Omega(\lambda,t),\lambda,t)d\lambda \,.$$

In order to complete the construction of the solution, we have to compute the integrals of $\mathfrak{G}_3$ and $\mathfrak{G}_2$, which appear in (4.10), in terms of the data. In order to achieve this, we will use the equations (4,5), (4,6) and (4,7), which can be written as follows:

(4.11) $e^{i\lambda x}\hat{u}(-\lambda,t) = e^{i\lambda x - \Omega(\lambda,t)}\hat{u}_0(-\lambda) + \mathfrak{G}_3(\lambda,x,t) - i\lambda \mathfrak{G}_2(\lambda,x,t)$

$$- \lambda^2 \mathfrak{G}_1(\lambda,x,t) + i\lambda^3 \mathfrak{G}_0(\lambda,x,t) + e^{i\lambda x - \Omega(\lambda,t)} \widetilde{\tilde{f}}(\Omega(\lambda,t),-\lambda,t),\text{ for } \mathrm{Im}(-\lambda) \leq 0\,,$$

(4.12) $e^{i\lambda x}\hat{u}(i\lambda,t) = e^{i\lambda x - \Omega(\lambda,t)}\hat{u}_0(i\lambda) + \mathfrak{G}_3(\lambda,x,t) - \lambda \mathfrak{G}_2(\lambda,x,t)$

$$+ \lambda^2 \mathfrak{G}_1(\lambda,x,t) - \lambda^3 \mathfrak{G}_0(\lambda,x,t) + e^{i\lambda x - \Omega(\lambda,t)} \widetilde{\tilde{f}}(\Omega(\lambda,t),i\lambda,t),\text{ for } \mathrm{Im}(i\lambda) \leq 0\,,$$

(4.13) $e^{i\lambda x}\hat{u}(-i\lambda,t) = e^{i\lambda x - \Omega(\lambda,t)}\hat{u}_0(-i\lambda) + \mathfrak{G}_3(\lambda,x,t) + \lambda \mathfrak{G}_2(\lambda,x,t)$

$$+ \lambda^2 \mathfrak{G}_1(\lambda,x,t) + \lambda^3 \mathfrak{G}_0(\lambda,x,t) + e^{i\lambda x - \Omega(\lambda,t)} \widetilde{\tilde{f}}(\Omega(\lambda,t),-i\lambda,t)\,,\text{ for } \lambda \in \mathbb{C} \text{ with } \mathrm{Im}(-i\lambda) \leq 0\,.$$

In order to determine the values of the integrals

$$\int_{\mathcal{L}_1} \mathfrak{G}_3(\lambda,x,t)d\lambda \text{ and } \int_{\mathcal{L}_1} i\lambda \mathfrak{G}_2(\lambda,x,t)d\lambda\,,$$

we solve the system of algebraic equations (4.11) and (4.13), for the quantities $\mathfrak{G}_3(\lambda,x,t)$ and $\mathfrak{G}_2(\lambda,x,t)$, with $\lambda \in \mathcal{L}_1$, observing that if $\lambda \in \mathcal{L}_1$ then $\mathrm{Im}(-\lambda) \leq 0$ and $\mathrm{Im}(-i\lambda) \leq 0$.

Similarly, in order to determine the values of the integrals

$$\int_{\mathcal{L}_2} \mathfrak{G}_3(\lambda,x,t)d\lambda \text{ and } \int_{\mathcal{L}_2} i\lambda \mathfrak{G}_2(\lambda,x,t)d\lambda\,,$$

we solve the system of algebraic equations (4.11) and (4.12), for the quantities $\mathfrak{G}_3(\lambda,x,t)$ and $\mathfrak{G}_2(\lambda,x,t)$, with $\lambda \in \mathcal{L}_2$, observing that if $\lambda \in \mathcal{L}_2$ then $\mathrm{Im}(-\lambda) \leq 0$ and $\mathrm{Im}(i\lambda) \leq 0$.

*Solving the system of* (4.11) & (4.13), *for* $\lambda \in \mathcal{L}_1$ Since the determinat of the coefficients of the unknown quantities is

$$\det\begin{bmatrix} 1 & -i\lambda \\ 1 & \lambda \end{bmatrix} = \lambda(1+i)\,,$$

we have, $\lambda \in \mathcal{L}_1$,

$$\mathfrak{G}_3(\lambda,x,t) = \frac{1}{1+i}\begin{bmatrix} e^{i\lambda x}\hat{u}(-\lambda,t) - e^{i\lambda x - \Omega(\lambda,t)}\hat{u}_0(-\lambda) + \lambda^2 \mathfrak{G}_1(\lambda,x,t) - i\lambda^3 \mathfrak{G}_0(\lambda,x,t) - e^{i\lambda x - \Omega(\lambda,t)} \widetilde{\tilde{f}}(\Omega(\lambda,t),-\lambda,t) & -i \\ e^{i\lambda x}\hat{u}(-i\lambda,t) - e^{i\lambda x - \Omega(\lambda,t)}\hat{u}_0(-i\lambda) - \lambda^2 \mathfrak{G}_1(\lambda,x,t) - \lambda^3 \mathfrak{G}_0(\lambda,x,t) - e^{i\lambda x - \Omega(\lambda,t)} \widetilde{\tilde{f}}(\Omega(\lambda,t),-i\lambda,t) & 1 \end{bmatrix}$$

and

$$i\lambda \mathcal{G}_2(\lambda,x,t) = \frac{i}{1+i}\begin{bmatrix} 1 & e^{i\lambda x}\hat{u}(-\lambda,t) - e^{i\lambda x - \Omega(\lambda,t)}\hat{u}_0(-\lambda) + \lambda^2 \mathcal{G}_1(\lambda,x,t) - i\lambda^3 \mathcal{G}_0(\lambda,x,t) - e^{i\lambda x - \Omega(\lambda,t)}\widetilde{\hat{f}}(\Omega(\lambda,t),-\lambda,t) \\ 1 & e^{i\lambda x}\hat{u}(-i\lambda,t) - e^{i\lambda x - \Omega(\lambda,t)}\hat{u}_0(-i\lambda) - \lambda^2 \mathcal{G}_1(\lambda,x,t) - \lambda^3 \mathcal{G}_0(\lambda,x,t) - e^{i\lambda x - \Omega(\lambda,t)}\widetilde{\hat{f}}(\Omega(\lambda,t),-i\lambda,t) \end{bmatrix}.$$

Thus, keeping in mind that

$$\int_{\mathcal{L}_1} e^{i\lambda x}\hat{u}(-\lambda,t)d\lambda = 0 \text{ and } \int_{\mathcal{L}_1} e^{i\lambda x}\hat{u}(-i\lambda,t)d\lambda = 0,$$

we obtain

$$(4.14) \quad \int_{\mathcal{L}_1} \mathcal{G}_3(\lambda,x,t)d\lambda = -\frac{1}{1+i}\int_{\mathcal{L}_1} e^{i\lambda x - \Omega(\lambda,t)}[\hat{u}_0(-\lambda) + i\hat{u}_0(-i\lambda)]d\lambda$$

$$+\frac{1}{1+i}\int_{\mathcal{L}_1}[(1-i)\lambda^2 \mathcal{G}_1(\lambda,x,t) - 2\lambda^3 \mathcal{G}_0(\lambda,x,t)]d\lambda - \frac{1}{1+i}\int_{\mathcal{L}_1} e^{i\lambda x - \Omega(\lambda,t)}[\widetilde{\hat{f}}(\Omega(\lambda,t),-\lambda,t) + i\widetilde{\hat{f}}(\Omega(\lambda,t),-i\lambda,t)]d\lambda$$

and

$$(4.15) \quad \int_{\mathcal{L}_1} i\lambda \mathcal{G}_2(\lambda,x,t)d\lambda = \frac{i}{1+i}\int_{\mathcal{L}_1} e^{i\lambda x - \Omega(\lambda,t)}[\hat{u}_0(-\lambda) - \hat{u}_0(-i\lambda)]d\lambda$$

$$-\frac{i}{1+i}\int_{\mathcal{L}_1}[2\lambda^2 \mathcal{G}_1(\lambda,x,t) + (1-i)\lambda^3 \mathcal{G}_0(\lambda,x,t)]d\lambda + \frac{i}{1+i}\int_{\mathcal{L}_1} e^{i\lambda x - \Omega(\lambda,t)}[\widetilde{\hat{f}}(\Omega(\lambda,t),-\lambda,t) - \widetilde{\hat{f}}(\Omega(\lambda,t),-i\lambda,t)]d\lambda.$$

*Solving the system of* (4.11)&(4.12), *for* $\lambda \in \mathcal{L}_2$ The determinat of the coefficients of the unknown quantities is

$$\det\begin{bmatrix} 1 & -i\lambda \\ 1 & -\lambda \end{bmatrix} = \lambda(-1+i),$$

and, for $\lambda \in \mathcal{L}_2$, we have

$$\mathcal{G}_3(\lambda,x,t) = \frac{1}{-1+i}\begin{bmatrix} e^{i\lambda x}\hat{u}(-\lambda,t) - e^{i\lambda x - \Omega(\lambda,t)}\hat{u}_0(-\lambda) + \lambda^2 \mathcal{G}_1(\lambda,x,t) - i\lambda^3 \mathcal{G}_0(\lambda,x,t) - e^{i\lambda x - \Omega(\lambda,t)}\widetilde{\hat{f}}(\Omega(\lambda,t),-\lambda,t) & -i \\ e^{i\lambda x}\hat{u}(i\lambda,t) - e^{i\lambda x - \Omega(\lambda,t)}\hat{u}_0(i\lambda) - \lambda^2 \mathcal{G}_1(\lambda,x,t) + \lambda^3 \mathcal{G}_0(\lambda,x,t) - e^{i\lambda x - \Omega(\lambda,t)}\widetilde{\hat{f}}(\Omega(\lambda,t),i\lambda,t) & -1 \end{bmatrix}$$

and

$$i\lambda \mathcal{G}_2(\lambda,x,t) = \frac{i}{-1+i}\begin{bmatrix} 1 & e^{i\lambda x}\hat{u}(-\lambda,t) - e^{i\lambda x - \Omega(\lambda,t)}\hat{u}_0(-\lambda) + \lambda^2 \mathcal{G}_1(\lambda,x,t) - i\lambda^3 \mathcal{G}_0(\lambda,x,t) - e^{i\lambda x - \Omega(\lambda,t)}\widetilde{\hat{f}}(\Omega(\lambda,t),-\lambda,t) \\ 1 & e^{i\lambda x}\hat{u}(i\lambda,t) - e^{i\lambda x - \Omega(\lambda,t)}\hat{u}_0(i\lambda) - \lambda^2 \mathcal{G}_1(\lambda,x,t) + \lambda^3 \mathcal{G}_0(\lambda,x,t) - e^{i\lambda x - \Omega(\lambda,t)}\widetilde{\hat{f}}(\Omega(\lambda,t),i\lambda,t) \end{bmatrix}.$$

Also,

$$\int_{\mathcal{L}_2} e^{i\lambda x}\hat{u}(-\lambda,t)d\lambda = 0 \text{ and } \int_{\mathcal{L}_2} e^{i\lambda x}\hat{u}(i\lambda,t)d\lambda = 0,$$

and, therefore,

$$(4.16) \quad \int_{\mathcal{L}_2} \mathcal{G}_3(\lambda,x,t)d\lambda = \frac{1}{-1+i}\int_{\mathcal{L}_2} e^{i\lambda x - \Omega(\lambda,t)}[\hat{u}_0(-\lambda) - i\hat{u}_0(i\lambda)]d\lambda$$

$$+\frac{1}{-1+i}\int_{\mathcal{L}_2}[-(1+i)\lambda^2 \mathcal{G}_1(\lambda,x,t) + 2i\lambda^3 \mathcal{G}_0(\lambda,x,t)]d\lambda + \frac{1}{-1+i}\int_{\mathcal{L}_2} e^{i\lambda x - \Omega(\lambda,t)}[\widetilde{\hat{f}}(\Omega(\lambda,t),-\lambda,t) - i\widetilde{\hat{f}}(\Omega(\lambda,t),i\lambda,t)]d\lambda$$

and

$$(4.17) \quad \int_{\mathcal{L}_2} i\lambda \mathcal{G}_2(\lambda,x,t)d\lambda = \frac{i}{-1+i}\int_{\mathcal{L}_2} e^{i\lambda x - \Omega(\lambda,t)}[\hat{u}_0(-\lambda) - \hat{u}_0(i\lambda)]d\lambda$$

$$-\frac{i}{-1+i}\int_{\mathcal{L}_2}[2\lambda^2\mathcal{G}_1(\lambda,x,t)-(1+i)\lambda^3\mathcal{G}_0(\lambda,x,t)]d\lambda + \frac{i}{-1+i}\int_{\mathcal{L}_2}e^{i\lambda x-\Omega(\lambda,t)}[\widetilde{\hat{f}}(\Omega(\lambda,t),-\lambda,t)-\widetilde{\hat{f}}(\Omega(\lambda,t),i\lambda,t)]d\lambda.$$

Finally, substituting (4.14), (4.15), 4.16) and (4.17), in (4.10), we obtain the formula for the solution of problem 3.

## 5. The operator $\partial_t + \alpha \partial_{xxxxt} - \beta(t)\partial_{xx}$

**Problem 4** Solve

(5.1)
$$\begin{cases} \partial_t u = -\alpha\partial_{xxxxt}u + \beta(t)\partial_{xx}u + f(x,t), \ (x,t)\in\mathbb{R}^+\times\mathbb{R}^+, \\ u(x,0)=u_0(x), \ x\in\mathbb{R}^+, \\ u(0,t)=g_0(t), \ t\in\mathbb{R}^+, \\ u_x(0,t)=g_1(t), \ t\in\mathbb{R}^+, \end{cases}$$

for $u(x,t)$.

***Derivation of the solution*** With notation as in the previous sections, taking Fourier transforms and integrating by parts, we see that the diferential equation in (5.1) gives, for $\lambda\in\mathbb{C}$ with $\mathrm{Im}\,\lambda\le 0$,

$$\frac{\partial}{\partial t}[\hat{u}(\lambda,t)]+\alpha\left\{\lambda^4\frac{\partial}{\partial t}[\hat{u}(\lambda,t)]-g'_3(t)-i\lambda g'_2(t)-(i\lambda)^2 g'_1(t)-(i\lambda)^3 g'_0(t)\right\}$$
$$+\beta(t)[\lambda^2\hat{u}(\lambda,t)+g_1(t)+i\lambda g_0(t)]=\hat{f}(\lambda,t),$$

$\Rightarrow$

$$(1+\alpha\lambda^4)\frac{\partial}{\partial t}[\hat{u}(\lambda,t)]+\beta(t)\lambda^2\hat{u}(\lambda,t)=\alpha[g'_3(t)+i\lambda g'_2(t)-\lambda^2 g'_1(t)-i\lambda^3 g'_0(t)]-\beta(t)[g_1(t)+i\lambda g_0(t)]+\hat{f}(\lambda,t)$$

$\Rightarrow$

$$\frac{\partial}{\partial t}[\hat{u}(\lambda,t)]+\frac{\beta(t)\lambda^2}{1+\alpha\lambda^4}\hat{u}(\lambda,t)=\frac{\alpha g'_3(t)}{1+\alpha\lambda^4}+i\lambda\frac{\alpha g'_2(t)}{1+\alpha\lambda^4}$$
$$-\frac{1}{1+\alpha\lambda^4}\{\alpha[\lambda^2 g'_1(t)+i\lambda^3 g'_0(t)]+\beta(t)[g_1(t)+i\lambda g_0(t)]\}+\frac{1}{1+\alpha\lambda^4}\hat{f}(\lambda,t)$$

$\Rightarrow$

(5.2)
$$\frac{\partial}{\partial t}[\hat{u}(\lambda,t)]+\frac{\beta(t)\lambda^2}{1+\alpha\lambda^4}\hat{u}(\lambda,t)=\frac{\alpha g'_3(t)}{1+\alpha\lambda^4}+i\lambda\frac{\alpha g'_2(t)}{1+\alpha\lambda^4}-\frac{G(\lambda,t)}{1+\alpha\lambda^4}+\frac{1}{1+\alpha\lambda^4}\hat{f}(\lambda,t),$$

where we have set

$$G(\lambda,t):=\alpha[\lambda^2 g'_1(t)+i\lambda^3 g'_0(t)]+\beta(t)[g_1(t)+i\lambda g_0(t)].$$

Also, (5.2) can be written as

(5.3)
$$\frac{\partial}{\partial t}[e^{\Omega(\lambda,t)}\hat{u}(\lambda,t)]=e^{\Omega(\lambda,t)}\frac{\alpha g'_3(t)}{1+\alpha\lambda^4}+i\lambda e^{\Omega(\lambda,t)}\frac{\alpha g'_2(t)}{1+\alpha\lambda^4}-e^{\Omega(\lambda,t)}\frac{G(\lambda,t)}{1+\alpha\lambda^4}+\frac{e^{\Omega(\lambda,t)}}{1+\alpha\lambda^4}\hat{f}(\lambda,t),$$

where

$$\omega(\lambda,t)=\frac{\beta(t)\lambda^2}{1+\alpha\lambda^4} \quad\text{and}\quad \Omega(\lambda,t)=\int_{\tau=0}^{t}\omega(\lambda,\tau)d\tau=\int_{\tau=0}^{t}\frac{\beta(\tau)\lambda^2}{1+\alpha\lambda^4}d\tau=\frac{\lambda^2}{1+\alpha\lambda^4}\int_{\tau=0}^{t}\beta(\tau)d\tau.$$

Integrating (5.3) we obtain

$$e^{\Omega(\lambda,t)}\hat{u}(\lambda,t)=\hat{u}_0(\lambda)+\frac{\alpha}{1+\alpha\lambda^4}\int_{\tau=0}^{t}e^{\Omega(\lambda,\tau)}g'_3(\tau)d\tau+\frac{\alpha i\lambda}{1+\alpha\lambda^4}\int_{\tau=0}^{t}e^{\Omega(\lambda,\tau)}g'_2(\tau)d\tau$$

$$-\frac{1}{1+\alpha\lambda^4}\int_{\tau=0}^{t}e^{\Omega(\lambda,\tau)}G(\lambda,\tau)d\tau+\int_{\tau=0}^{t}\frac{e^{\Omega(\lambda,t)}}{1+\alpha\lambda^4}\hat{f}(\lambda,t)d\tau,$$

or, equivalently,

(5.4) $\quad \hat{u}(\lambda,t) = e^{-\Omega(\lambda,t)}\hat{u}_0(\lambda) + \dfrac{\alpha}{1+\alpha\lambda^4}e^{-\Omega(\lambda,t)}(g'_3\,\widetilde{\,})(\Omega,t) + \dfrac{\alpha i\lambda}{1+\alpha\lambda^4}e^{-\Omega(\lambda,t)}(g'_2\,\widetilde{\,})(\Omega,t)$

$$-\frac{e^{-\Omega(\lambda,t)}}{1+\alpha\lambda^4}\widetilde{G}(\Omega,\lambda,t)+\frac{e^{-\Omega(\lambda,t)}}{1+\alpha\lambda^4}\widetilde{\hat{f}}(\Omega,\lambda,t),$$

where

$$(g'_3\,\widetilde{\,})(\Omega,t) = (g'_3\,\widetilde{\,})(\Omega(\lambda,t),t) := \int_{\tau=0}^{t}e^{\Omega(\lambda,\tau)}g'_3(\tau)d\tau,\quad (g'_2\,\widetilde{\,})(\Omega,t) = (g'_2\,\widetilde{\,})(\Omega(\lambda,t),t) := \int_{\tau=0}^{t}e^{\Omega(\lambda,\tau)}g'_2(\tau)d\tau,$$

$$\widetilde{G}(\Omega(\lambda,t),\lambda,t) = \widetilde{G}(\Omega,\lambda,t) := \int_{\tau=0}^{t}e^{\Omega(\lambda,\tau)}G(\lambda,\tau)d\tau \text{ and } \widetilde{\hat{f}}(\Omega(\lambda,t),\lambda,t) = \widetilde{\hat{f}}(\Omega,\lambda,t) := \int_{\tau=0}^{t}e^{\Omega(\lambda,\tau)}\hat{f}(\lambda,t)d\tau.$$

Multiplying (5.4) by $e^{i\lambda x}$ and integrating, we obtain

(5.5) $\quad 2\pi u(x,t) = \displaystyle\int_{-\infty}^{\infty}e^{i\lambda x}\hat{u}(\lambda,t)d\lambda = \int_{-\infty}^{\infty}e^{i\lambda x-\Omega(\lambda,t)}\hat{u}_0(\lambda)d\lambda$

$$+\alpha\int_{-\infty}^{\infty}\frac{e^{i\lambda x-\Omega(\lambda,t)}}{1+\alpha\lambda^4}(g'_3\,\widetilde{\,})(\Omega,t)d\lambda+\alpha i\int_{-\infty}^{\infty}\frac{\lambda e^{i\lambda x-\Omega(\lambda,t)}}{1+\alpha\lambda^4}(g'_2\,\widetilde{\,})(\Omega,t)d\lambda$$

$$-\int_{-\infty}^{\infty}\frac{e^{i\lambda x-\Omega(\lambda,t)}}{1+\alpha\lambda^4}\widetilde{G}(\Omega,\lambda,t)d\lambda+\int_{-\infty}^{\infty}\frac{e^{i\lambda x-\Omega(\lambda,t)}}{1+\alpha\lambda^4}\widetilde{\hat{f}}(\Omega,\lambda,t)d\lambda.$$

Next, observing that

$$1+\alpha\lambda^4 = 0 \;\Rightarrow\; \lambda = \pm e^{i\pi/4}/\sqrt[4]{\alpha} \text{ or } \lambda = \pm e^{i3\pi/4}/\sqrt[4]{\alpha},$$

we consider the numbers

$$\lambda_1 := e^{i\pi/4}/\sqrt[4]{\alpha} \text{ and } \lambda_2 := e^{i3\pi/4}/\sqrt[4]{\alpha}.$$

Deforming the contours of the integrals in (5.5), which involve the quantities $(g'_3\,\widetilde{\,})(\Omega,t)$, $(g'_2\,\widetilde{\,})(\Omega,t)$ and $\widetilde{G}(\Omega,\lambda,t)$, we obtain

(5.6) $\quad 2\pi u(x,t) = \displaystyle\int_{-\infty}^{\infty}e^{i\lambda x-\Omega(\lambda,t)}\hat{u}_0(\lambda)d\lambda$

$$+\alpha\int_{C_1+C_2}\frac{e^{i\lambda x-\Omega(\lambda,t)}}{1+\alpha\lambda^4}(g'_3\,\widetilde{\,})(\Omega,t)d\lambda+\alpha i\int_{C_1+C_2}\frac{\lambda e^{i\lambda x-\Omega(\lambda,t)}}{1+\alpha\lambda^4}(g'_2\,\widetilde{\,})(\Omega,t)d\lambda$$

$$-\int_{C_1+C_2}\frac{e^{i\lambda x-\Omega(\lambda,t)}}{1+\alpha\lambda^4}\widetilde{G}(\Omega,\lambda,t)d\lambda+\int_{-\infty}^{\infty}\frac{e^{i\lambda x-\Omega(\lambda,t)}}{1+\alpha\lambda^4}\widetilde{\hat{f}}(\Omega,\lambda,t)d\lambda,$$

where the contours $C_1$ and $C_2$ are simple closed curves, around the points $\lambda_1$ and $\lambda_2$, as in fig. 5.

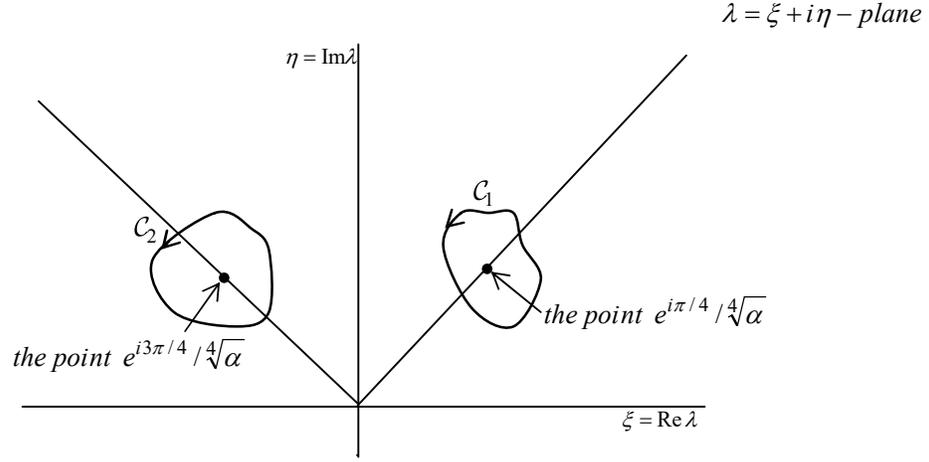

Fig. 5 Choices of the contours $C_1$ and $C_2$.

In order to complete the construction of the solution, we have to compute the integrals in (5.6), which contain $(g'_3\widetilde{\;})(\Omega,t)$ and $(g'_2\widetilde{\;})(\Omega,t)$, in terms of the data. In order to achieve this, first we solve the equation

$$\frac{[\sigma(\lambda)]^2}{1+\alpha[\sigma(\lambda)]^4} = \frac{\lambda^2}{1+\alpha\lambda^4}$$

for $\sigma(\lambda)$, and we find

$$\sigma_0(\lambda)=\lambda,\ \sigma_1(\lambda)=-\lambda,\ \sigma_2(\lambda)=1/(\lambda\sqrt{\alpha}),\ \sigma_3(\lambda)=-1/(\lambda\sqrt{\alpha}).$$

For these functions we have

$$\omega(\sigma_j(\lambda),t)=\omega(\lambda,t)\ \text{and}\ \Omega(\sigma_j(\lambda),t)=\Omega(\lambda,t),\ j=0,1,2,3.$$

Also, we observe that

$$\sigma_1(\lambda_1)=-\lambda_1,\ \sigma_1(\lambda_2)=-\lambda_2,\ \sigma_2(\lambda_1)=\frac{e^{-i\pi/4}}{\sqrt{\alpha}}\ \text{and}\ \sigma_2(\lambda_2)=\frac{e^{-i3\pi/4}}{\sqrt{\alpha}}.$$

Thus,

$$\operatorname{Im}\sigma_1(\lambda_1)<0,\ \sigma_1(\lambda_2)<0,\ \sigma_2(\lambda_1)<0,\ \text{and}\ \sigma_2(\lambda_2)<0.$$

and, therefore, (5.4) imply

(5.7)  $\hat{u}(-\lambda,t)=e^{-\Omega(\lambda,t)}\hat{u}_0(-\lambda)+\dfrac{\alpha}{1+\alpha\lambda^4}e^{-\Omega(\lambda,t)}(g'_3\widetilde{\;})(\Omega,t)-\dfrac{\alpha i\lambda}{1+\alpha\lambda^4}e^{-\Omega(\lambda,t)}(g'_2\widetilde{\;})(\Omega,t)$

$$-\frac{e^{-\Omega(\lambda,t)}}{1+\alpha\lambda^4}\widetilde{G}(\Omega,-\lambda,t)+\frac{e^{-\Omega(\lambda,t)}}{1+\alpha\lambda^4}\widetilde{\tilde{f}}(\Omega,-\lambda,t),\ \text{for}\ \lambda\in C(\lambda_1)\cup C(\lambda_2),$$

and

(5.8)  $\hat{u}(1/(\lambda\sqrt{\alpha}),t)=e^{-\Omega(\lambda,t)}\hat{u}_0(1/(\lambda\sqrt{\alpha}))+\dfrac{\alpha^2\lambda^4}{1+\alpha\lambda^4}e^{-\Omega(\lambda,t)}(g'_3\widetilde{\;})(\Omega,t)+\dfrac{\alpha^2 i\lambda^4}{1+\alpha\lambda^4}e^{-\Omega(\lambda,t)}(g'_2\widetilde{\;})(\Omega,t)$

$$-\frac{\alpha\lambda^4}{1+\alpha\lambda^4}e^{-\Omega(\lambda,t)}\widetilde{G}(\Omega,1/(\lambda\sqrt{\alpha}),t)+\frac{\alpha\lambda^4}{1+\alpha\lambda^4}e^{-\Omega(\lambda,t)}\widetilde{\tilde{f}}(\Omega,1/(\lambda\sqrt{\alpha}),t),\ \text{for}\ \lambda\in C(\lambda_1)\cup C(\lambda_2).$$

(*Note* The function $\sigma_3(\lambda)=-1/(\lambda\sqrt{\alpha})$ is of no use for our purposes, since $\sigma_3(\lambda_1)=-e^{-i\pi/4}/\sqrt{\alpha}$ and $\sigma_3(\lambda_2)=-e^{-i3\pi/4}/\sqrt{\alpha}$, hence $\operatorname{Im}\sigma_3(\lambda_1)>0$ and $\operatorname{Im}\sigma_3(\lambda_2)>0$.)

Next, we solve the system of algebraic equations (5.7) & (5.8), for the unkown quantities

$$\frac{\alpha}{1+\alpha\lambda^4}e^{-\Omega(\lambda,t)}(g'_3\,\widetilde{\,})(\Omega,t) \text{ and } \frac{\alpha}{1+\alpha\lambda^4}e^{-\Omega(\lambda,t)}(g'_2\,\widetilde{\,})(\Omega,t), \text{ with } \lambda \in \mathcal{C}(\lambda_1)\cup\mathcal{C}(\lambda_2).$$

Computing the determinant of the coefficients of the unknowns, namely

$$\det\begin{bmatrix} 1 & -i\lambda \\ \alpha\lambda^4 & \alpha i\lambda^4 \end{bmatrix} = \alpha i\lambda^4(1+\lambda),$$

we find that, for $\lambda \in \mathcal{C}(\lambda_1)\cup\mathcal{C}(\lambda_2)$,

(5.9) $\dfrac{\alpha}{1+\alpha\lambda^4}e^{-\Omega(\lambda,t)}(g'_3\,\widetilde{\,})(\Omega,t)$

$$= \frac{1}{\alpha i\lambda^4(1+\lambda)}\det\begin{bmatrix} \hat{u}(-\lambda,t)-e^{-\Omega(\lambda,t)}\hat{u}_0(-\lambda)+\dfrac{e^{-\Omega(\lambda,t)}}{1+\alpha\lambda^4}\widetilde{G}(\Omega,-\lambda,t)-\dfrac{e^{-\Omega(\lambda,t)}}{1+\alpha\lambda^4}\widetilde{\tilde{f}}(\Omega,-\lambda,t) & -i\lambda \\ \hat{u}(1/(\lambda\sqrt{\alpha}),t)-e^{-\Omega(\lambda,t)}\hat{u}_0(1/(\lambda\sqrt{\alpha}))+\dfrac{\alpha\lambda^4 e^{-\Omega(\lambda,t)}}{1+\alpha\lambda^4}[\widetilde{G}(\Omega,1/(\lambda\sqrt{\alpha}),t)-\widetilde{\tilde{f}}(\Omega,1/(\lambda\sqrt{\alpha}),t)] & \alpha i\lambda^4 \end{bmatrix}$$

and

(5.10) $\dfrac{\alpha}{1+\alpha\lambda^4}e^{-\Omega(\lambda,t)}(g'_2\,\widetilde{\,})(\Omega,t)$

$$= \frac{1}{\alpha i\lambda^4(1+\lambda)}\det\begin{bmatrix} 1 & \hat{u}(-\lambda,t)-e^{-\Omega(\lambda,t)}\hat{u}_0(-\lambda)+\dfrac{e^{-\Omega(\lambda,t)}}{1+\alpha\lambda^4}\widetilde{G}(\Omega,-\lambda,t)-\dfrac{e^{-\Omega(\lambda,t)}}{1+\alpha\lambda^4}\widetilde{\tilde{f}}(\Omega,-\lambda,t) \\ \alpha i\lambda^4 & \hat{u}(1/(\lambda\sqrt{\alpha}),t)-e^{-\Omega(\lambda,t)}\hat{u}_0(1/(\lambda\sqrt{\alpha}))+\dfrac{\alpha\lambda^4 e^{-\Omega(\lambda,t)}}{1+\alpha\lambda^4}[\widetilde{G}(\Omega,1/(\lambda\sqrt{\alpha}),t)-\widetilde{\tilde{f}}(\Omega,1/(\lambda\sqrt{\alpha}),t)] \end{bmatrix}.$$

Substituting (5.9) and (5.10) in (5.6), and taking into consideration that

$$\int_{C_1+C_2}\frac{1}{\lambda+1}e^{i\lambda x}\hat{u}(-\lambda,t)d\lambda = 0, \quad \int_{C_1+C_2}\frac{\lambda}{\lambda+1}e^{i\lambda x}\hat{u}(-\lambda,t)d\lambda = 0, \quad \int_{C_1+C_2}e^{i\lambda x}\frac{1}{\lambda^3(1+\lambda)}\hat{u}(1/(\lambda\sqrt{\alpha}),t)d\lambda = 0,$$

we obtain the formula for the solution of problem 4:

(5.11)  $2\pi u(x,t) = \displaystyle\int_{-\infty}^{\infty}e^{i\lambda x-\Omega(\lambda,t)}\hat{u}_0(\lambda)d\lambda - \int_{C_1+C_2}\frac{1-i\lambda}{1+\lambda}e^{i\lambda x-\Omega(\lambda,t)}\hat{u}_0(-\lambda)d\lambda -\frac{2}{\alpha}\int_{C_1+C_2}\frac{e^{i\lambda x-\Omega(\lambda,t)}}{\lambda^3(1+\lambda)}\hat{u}_0(1/(\lambda\sqrt{\alpha}))d\lambda$

$\displaystyle -\int_{C_1+C_2}\frac{e^{i\lambda x-\Omega(\lambda,t)}}{1+\alpha\lambda^4}\widetilde{G}(\Omega,\lambda,t)d\lambda + \int_{C_1+C_2}\frac{1-i\lambda}{1+\lambda}\frac{e^{i\lambda x-\Omega(\lambda,t)}}{1+\alpha\lambda^4}\widetilde{G}(\Omega,-\lambda,t)d\lambda + \int_{C_1+C_2}\frac{2\lambda}{1+\lambda}\frac{e^{i\lambda x-\Omega(\lambda,t)}}{1+\alpha\lambda^4}\widetilde{G}(\Omega,1/(\lambda\sqrt{\alpha}),t)d\lambda$

$\displaystyle +\int_{-\infty}^{\infty}\frac{e^{i\lambda x-\Omega(\lambda,t)}}{1+\alpha\lambda^4}\widetilde{\tilde{f}}(\Omega,\lambda,t)d\lambda - \int_{C_1+C_2}\frac{1-i\lambda}{1+\lambda}\frac{e^{i\lambda x-\Omega(\lambda,t)}}{1+\alpha\lambda^4}\widetilde{\tilde{f}}(\Omega,-\lambda,t)d\lambda - \int_{C_1+C_2}\frac{2\lambda}{1+\lambda}\frac{e^{i\lambda x-\Omega(\lambda,t)}}{1+\alpha\lambda^4}\widetilde{\tilde{f}}(\Omega,1/(\lambda\sqrt{\alpha}),t)d\lambda d\lambda.$

## 6. The operator $\partial_t + \alpha\partial_{xxxxt} - \beta\partial_{xxt} + \gamma(t)\partial_{xxxx}$

**Problem 5** Solve

(6.1) $\begin{cases} \partial_t u = -\alpha\partial_{xxxxt}u + \beta\partial_{xxt}u - \gamma(t)\partial_{xxxx}u + f(x,t), \ (x,t)\in\mathbb{R}^+\times\mathbb{R}^+, \\ u(x,0) = u_0(x), \ x\in\mathbb{R}^+, \\ u(0,t) = g_0(t), \ t\in\mathbb{R}^+, \\ u_x(0,t) = g_1(t), \ t\in\mathbb{R}^+, \end{cases}$

*for* $u(x,t)$.

***Derivation of the solution*** As in the previous sections, we see that, for $\lambda \in \mathbb{C}$ with $\operatorname{Im}\lambda \leq 0$, the differential equation in (6.1) leads to

$$(1+\beta\lambda^2+\alpha\lambda^4)\frac{\partial}{\partial t}[\hat{u}(\lambda,t)]+\gamma(t)\lambda^4\hat{u}(\lambda,t) = \alpha[g'_3(t)+i\lambda g'_2(t)-\lambda^2 g'_1(t)-i\lambda^3 g'_0(t)]$$

$$-\beta[g'_1(t)+i\lambda g'_0(t)]+\gamma(t)[g_3(t)+i\lambda g_2(t)-\lambda^2 g_1(t)-i\lambda^3 g_0(t)]+\hat{f}(\lambda,t)$$

$\Rightarrow$

$$\frac{\partial}{\partial t}[\hat{u}(\lambda,t)]+\frac{\gamma(t)\lambda^4}{1+\beta\lambda^2+\alpha\lambda^4}\hat{u}(\lambda,t) = \frac{\alpha}{1+\beta\lambda^2+\alpha\lambda^4}[g'_3(t)+i\lambda g'_2(t)-\lambda^2 g'_1(t)-i\lambda^3 g'_0(t)]$$

$$-\frac{\beta}{1+\beta\lambda^2+\alpha\lambda^4}[g'_1(t)+i\lambda g'_0(t)]+\frac{\gamma(t)}{1+\beta\lambda^2+\alpha\lambda^4}[g_3(t)+i\lambda g_2(t)-\lambda^2 g_1(t)-i\lambda^3 g_0(t)]+\frac{\hat{f}(\lambda,t)}{1+\beta\lambda^2+\alpha\lambda^4}$$

$\Rightarrow$

$$(6.2) \quad \frac{\partial}{\partial t}[\hat{u}(\lambda,t)]+\frac{\gamma(t)\lambda^4}{1+\beta\lambda^2+\alpha\lambda^4}\hat{u}(\lambda,t)$$

$$= \frac{1}{1+\beta\lambda^2+\alpha\lambda^4}[\alpha g'_3(t)+\gamma(t)g_3(t)]+\frac{i\lambda}{1+\beta\lambda^2+\alpha\lambda^4}[\alpha g'_2(t)+\gamma(t)g_2(t)]$$

$$+\frac{G(\lambda,t)}{1+\beta\lambda^2+\alpha\lambda^4}+\frac{\hat{f}(\lambda,t)}{1+\beta\lambda^2+\alpha\lambda^4},$$

where we have set

$$G(\lambda,t):=-\alpha[\lambda^2 g'_1(t)+i\lambda^3 g'_0(t)]-\beta[g'_1(t)+i\lambda g'_0(t)]-\gamma(t)[\lambda^2 g_1(t)+i\lambda^3 g_0(t)]$$

$$=-(\alpha\lambda^2+\beta)g'_1(t)-i\lambda(\alpha\lambda^2+\beta)g'_0(t)-\gamma(t)\lambda^2[g_1(t)+i\lambda g_0(t)].$$

Furthermore, setting

$$\omega(\lambda,t)=\frac{\gamma(t)\lambda^4}{1+\beta\lambda^2+\alpha\lambda^4}, \quad \Omega(\lambda,t)=\int_{\tau=0}^{t}\omega(\lambda,\tau)d\tau=\int_{\tau=0}^{t}\frac{\gamma(\tau)\lambda^4}{1+\beta\lambda^2+\alpha\lambda^4}d\tau=\frac{\lambda^4}{1+\beta\lambda^2+\alpha\lambda^4}\int_{\tau=0}^{t}\gamma(\tau)d\tau,$$

and

$$G_3(t)=\alpha g'_3(t)+\gamma(t)g_3(t) \text{ and } G_2(t)=\alpha g'_2(t)+\gamma(t)g_2(t),$$

we write (6.2) in the following way:

$$\frac{\partial}{\partial t}[e^{\Omega(\lambda,t)}\hat{u}(\lambda,t)]=\frac{e^{\Omega(\lambda,t)}G_3(t)}{1+\beta\lambda^2+\alpha\lambda^4}+\frac{i\lambda e^{\Omega(\lambda,t)}G_2(t)}{1+\beta\lambda^2+\alpha\lambda^4}-\frac{e^{\Omega(\lambda,t)}G(\lambda,t)}{1+\beta\lambda^2+\alpha\lambda^4}+\frac{e^{\Omega(\lambda,t)}\hat{f}(\lambda,t)}{1+\beta\lambda^2+\alpha\lambda^4}.$$

Integrating the above equation, we find

$$(6.3) \quad e^{\Omega(\lambda,t)}\hat{u}(\lambda,t)=\hat{u}_0(\lambda)+\int_{\tau=0}^{t}\frac{e^{\Omega(\lambda,\tau)}G_3(\tau)}{1+\beta\lambda^2+\alpha\lambda^4}d\tau+\int_{\tau=0}^{t}\frac{i\lambda e^{\Omega(\lambda,\tau)}G_2(\tau)}{1+\beta\lambda^2+\alpha\lambda^4}d\tau$$

$$+\int_{\tau=0}^{t}\frac{e^{\Omega(\lambda,\tau)}G(\lambda,\tau)}{1+\beta\lambda^2+\alpha\lambda^4}d\tau+\int_{\tau=0}^{t}\frac{e^{\Omega(\lambda,\tau)}\hat{f}(\lambda,\tau)}{1+\beta\lambda^2+\alpha\lambda^4}d\tau.$$

Defining

$$\widetilde{G}_3(\Omega(\lambda,t),t)=\int_{\tau=0}^{t}e^{\Omega(\lambda,\tau)}G_3(\tau)d\tau, \quad \widetilde{G}_2(\Omega(\lambda,t),t)=\int_{\tau=0}^{t}e^{\Omega(\lambda,\tau)}G_2(\tau)d\tau,$$

$$\widetilde{G}(\Omega(\lambda,t),\lambda,t) = \int_{\tau=0}^{t} e^{\Omega(\lambda,\tau)} G(\lambda,\tau) d\tau \text{ and } \widetilde{\hat{f}}(\Omega(\lambda,t),\lambda,t) = \int_{\tau=0}^{t} e^{\Omega(\lambda,\tau)} \hat{f}(\lambda,\tau) d\tau,$$

we write (6.3) as follows:

(6.4) $\quad \hat{u}(\lambda,t) = e^{-\Omega(\lambda,t)} \hat{u}_0(\lambda) + \dfrac{e^{-\Omega(\lambda,t)}}{1+\beta\lambda^2 + \alpha\lambda^4} \widetilde{G}_3(\Omega(\lambda,t),t) + \dfrac{i\lambda e^{-\Omega(\lambda,t)}}{1+\beta\lambda^2 + \alpha\lambda^4} \widetilde{G}_2(\Omega(\lambda,t),t)$

$$+ \frac{e^{-\Omega(\lambda,t)}}{1+\beta\lambda^2 + \alpha\lambda^4} \widetilde{G}(\Omega(\lambda,t),\lambda,t) + \frac{e^{-\Omega(\lambda,t)}}{1+\beta\lambda^2 + \alpha\lambda^4} \widetilde{\hat{f}}(\Omega(\lambda,t),\lambda,t), \text{ for } \lambda \in \mathbb{C} \text{ with } \mathrm{Im}\,\lambda \leq 0.$$

Now (6.4) gives

(6.5) $\quad 2\pi u(x,t) = \displaystyle\int_{-\infty}^{\infty} e^{i\lambda x} \hat{u}(\lambda,t) d\lambda$

$$= \int_{-\infty}^{\infty} e^{i\lambda x - \Omega(\lambda,t)} \hat{u}_0(\lambda) d\lambda + \int_{-\infty}^{\infty} \frac{e^{i\lambda x - \Omega(\lambda,t)}}{1+\beta\lambda^2 + \alpha\lambda^4} \widetilde{G}_3(\Omega(\lambda,t),t) d\lambda + \int_{-\infty}^{\infty} \frac{i\lambda e^{i\lambda x - \Omega(\lambda,t)}}{1+\beta\lambda^2 + \alpha\lambda^4} \widetilde{G}_2(\Omega(\lambda,t),t) d\lambda$$

$$+ \int_{-\infty}^{\infty} \frac{e^{i\lambda x - \Omega(\lambda,t)}}{1+\beta\lambda^2 + \alpha\lambda^4} \widetilde{G}(\Omega(\lambda,t),\lambda,t) d\lambda + \int_{-\infty}^{\infty} \frac{e^{i\lambda x - \Omega(\lambda,t)}}{1+\beta\lambda^2 + \alpha\lambda^4} \widetilde{\hat{f}}(\Omega(\lambda,t),\lambda,t) d\lambda.$$

We continue with the construction of the solution of problem (6.1), assuming $4\alpha > \beta^2$.

**The case** $4\alpha > \beta^2$  In this case the zeros of the polynomial $\Pi(\lambda) := 1+\beta\lambda^2 + \alpha\lambda^4$ are the numbers $\lambda_1$, $\lambda_2$, $-\lambda_1$, $-\lambda_2$, where

$$\lambda_1 = \sqrt{\frac{-\beta + i\sqrt{4\alpha - \beta^2}}{2\alpha}} \text{ and } \lambda_2 = \sqrt{\frac{-\beta - i\sqrt{4\alpha - \beta^2}}{2\alpha}}.$$

(The above square roots are chosen so that $\lambda_1$ lies in the first quarter of the plane and $\lambda_2$ lies in the second one.)

Next we see that any of the functions

$$\sigma(\lambda) = \pm\lambda \text{ or } \pm\sqrt{\frac{-\lambda^2}{1+\beta\lambda^2}},$$

satisfy the equation

$$\frac{\gamma(t)[\sigma(\lambda)]^4}{1+\beta[\sigma(\lambda)]^2 + \alpha[\sigma(\lambda)]^4} = \frac{\gamma(t)\lambda^4}{1+\beta\lambda^2 + \alpha\lambda^4}, \text{ i.e., } \omega(\sigma(\lambda),t) = \omega(\lambda,t).$$

Also,

$$\Omega(\sigma(\lambda),t) = \Omega(\lambda,t).$$

Now, we specify the square root and we define the function

$$\sigma_1(\lambda) := \sqrt{\frac{-\lambda^2}{1+\beta\lambda^2}}, \text{ to be analytic for } \lambda \text{ in an open neighborhood of } \lambda_1,$$

in such a way that

$$\sigma_1(\lambda_1) = \sqrt{\frac{-\lambda_1^2}{1+\beta\lambda_1^2}} = \sqrt{\frac{-\lambda_1^2}{-\alpha\lambda_1^4}} = \frac{1}{\lambda_1\sqrt{\alpha}}.$$

Similarly, we choose the square root and we define the function

$$\sigma_2(\lambda) := \sqrt{\frac{-\lambda^2}{1+\beta\lambda^2}}\,,\ \text{to be analytic for } \lambda \text{ in an open neighborhood of } \lambda_2,$$

such that

$$\sigma_2(\lambda_2) = \sqrt{\frac{-\lambda_2^2}{1+\beta\lambda_2^2}} = \sqrt{\frac{-\lambda_2^2}{-\alpha\lambda_2^4}} = \frac{1}{\lambda_2\sqrt{\alpha}}.$$

Let us note also that

(6.6) $$\Pi(\sigma_1(\lambda)) = \Pi(\sigma_2(\lambda)) = \frac{1+\beta\lambda^2+\alpha\lambda^4}{(1+\beta\lambda^2)^2} = \frac{\Pi(\lambda)}{(1+\beta\lambda^2)^2}$$

In particular, each of the quotients

$$\frac{\Pi(\sigma_1(\lambda))}{\Pi(\lambda)} \ \text{and} \ \frac{\Pi(\sigma_2(\lambda))}{\Pi(\lambda)}$$

is an analytic function for $\lambda$ in an open neighborhood of $\lambda_1$ and $\lambda_2$, respectively.

Considering contours $C_1$ and $C_2$, as the ones depicted in fig 6, we write (6.5) as follows:

(6.7) $$2\pi u(x,t) = \int_{-\infty}^{\infty} e^{i\lambda x - \Omega(\lambda,t)}\hat{u}_0(\lambda)d\lambda + \int_{C_1} \frac{e^{i\lambda x - \Omega(\lambda,t)}}{\Pi(\lambda)}\widetilde{G}_3(\Omega(\lambda,t),t)d\lambda + \int_{C_2} \frac{e^{i\lambda x - \Omega(\lambda,t)}}{\Pi(\lambda)}\widetilde{G}_3(\Omega(\lambda,t),t)d\lambda$$

$$+ \int_{C_1} \frac{i\lambda e^{i\lambda x - \Omega(\lambda,t)}}{\Pi(\lambda)}\widetilde{G}_2(\Omega(\lambda,t),t)d\lambda + \int_{C_2} \frac{i\lambda e^{i\lambda x - \Omega(\lambda,t)}}{\Pi(\lambda)}\widetilde{G}_2(\Omega(\lambda,t),t)d\lambda$$

$$+ \int_{C_1} \frac{e^{i\lambda x - \Omega(\lambda,t)}}{\Pi(\lambda)}\widetilde{G}(\Omega(\lambda,t),\lambda,t)d\lambda + \int_{C_2} \frac{e^{i\lambda x - \Omega(\lambda,t)}}{\Pi(\lambda)}\widetilde{G}(\Omega(\lambda,t),\lambda,t)d\lambda + \int_{-\infty}^{\infty} \frac{e^{i\lambda x - \Omega(\lambda,t)}}{\Pi(\lambda)}\widetilde{\hat{f}}(\Omega(\lambda,t),\lambda,t)d\lambda.$$

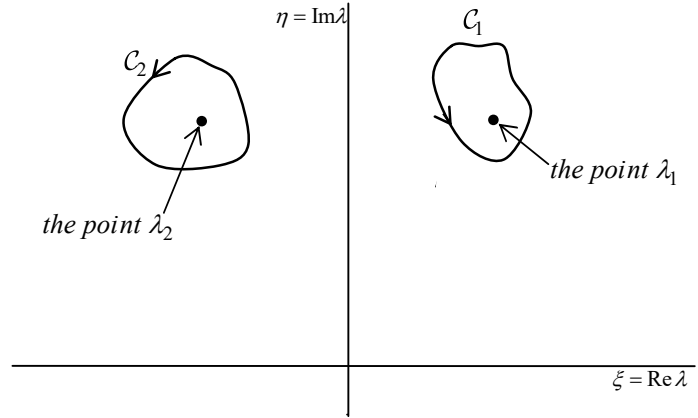

Fig. 6 Choices of the contours $C_1$ and $C_2$.

Next, from (6.4) we obtain

(6.8) $$\hat{u}(-\lambda,t) = e^{-\Omega(\lambda,t)}\hat{u}_0(-\lambda) + \frac{e^{-\Omega(\lambda,t)}}{\Pi(\lambda)}\widetilde{G}_3(\Omega(\lambda,t),t) - \frac{i\lambda e^{-\Omega(\lambda,t)}}{\Pi(\lambda)}\widetilde{G}_2(\Omega(\lambda,t),t)$$

$$+ \frac{e^{-\Omega(\lambda,t)}}{\Pi(\lambda)}\widetilde{G}(\Omega(\lambda,t),-\lambda,t) + \frac{e^{-\Omega(\lambda,t)}}{\Pi(\lambda)}\widetilde{\hat{f}}(\Omega(\lambda,t),-\lambda,t), \text{ for } \lambda \in C_1 \cup C_2,$$

$$(6.9) \quad \hat{u}(\sigma_1(\lambda),t) = e^{-\Omega(\lambda,t)}\hat{u}_0(\sigma_1(\lambda)) + \frac{e^{-\Omega(\lambda,t)}}{\Pi(\sigma_1(\lambda))}\widetilde{G}_3(\Omega(\lambda,t),t) + \frac{i\sigma_1(\lambda)e^{-\Omega(\lambda,t)}}{\Pi(\sigma_1(\lambda))}\widetilde{G}_2(\Omega(\lambda,t),t)$$

$$+ \frac{e^{-\Omega(\lambda,t)}}{\Pi(\sigma_1(\lambda))}\widetilde{G}(\Omega(\lambda,t),\sigma_1(\lambda),t) + \frac{e^{-\Omega(\lambda,t)}}{\Pi(\sigma_1(\lambda))}\widetilde{\tilde{f}}(\Omega(\lambda,t),\sigma_1(\lambda),t), \text{ for } \lambda \in \mathcal{C}_1,$$

$$(6.10) \quad \hat{u}(\sigma_2(\lambda),t) = e^{-\Omega(\lambda,t)}\hat{u}_0(\sigma_2(\lambda)) + \frac{e^{-\Omega(\lambda,t)}}{\Pi(\sigma_2(\lambda))}\widetilde{G}_3(\Omega(\lambda,t),t) + \frac{i\sigma_2(\lambda)e^{-\Omega(\lambda,t)}}{\Pi(\sigma_2(\lambda))}\widetilde{G}_2(\Omega(\lambda,t),t)$$

$$+ \frac{e^{-\Omega(\lambda,t)}}{\Pi(\sigma_2(\lambda))}\widetilde{G}(\Omega(\lambda,t),\sigma_2(\lambda),t) + \frac{e^{-\Omega(\lambda,t)}}{\Pi(\sigma_2(\lambda))}\widetilde{\tilde{f}}(\Omega(\lambda,t),\sigma_2(\lambda),t), \text{ for } \lambda \in \mathcal{C}_2.$$

From the equations (6.8), (6.9) and (6.10), we will determine the values of the quantities

$$(6.11) \quad e^{-\Omega(\lambda,t)}\widetilde{G}_3(\Omega(\lambda,t),t) \text{ and } e^{-\Omega(\lambda,t)}\widetilde{G}_2(\Omega(\lambda,t),t), \text{ for } \lambda \in \mathcal{C}_1 \text{ and for } \lambda \in \mathcal{C}_2.$$

*Solving the system of equations* (6.8) & (6.9) Working for $\lambda \in \mathcal{C}_1$, we compute the determinant of the unknowns (6.11):

$$\det\begin{bmatrix} \dfrac{1}{\Pi(\lambda)} & -\dfrac{i\lambda}{\Pi(\lambda)} \\ \dfrac{1}{\Pi(\sigma_1(\lambda))} & \dfrac{i\sigma_1(\lambda)}{\Pi(\sigma_1(\lambda))} \end{bmatrix} = \frac{i[\lambda+\sigma_1(\lambda)]}{\Pi(\lambda)\Pi(\sigma_1(\lambda))}.$$

Thus,

$$(6.12) \quad e^{-\Omega(\lambda,t)}\widetilde{G}_3(\Omega(\lambda,t),t) = \frac{\Pi(\lambda)\Pi(\sigma_1(\lambda))}{i[\lambda+\sigma_1(\lambda)]} \times$$

$$\times \det\begin{bmatrix} \hat{u}(-\lambda,t) - e^{-\Omega(\lambda,t)}\hat{u}_0(-\lambda) - \dfrac{e^{-\Omega(\lambda,t)}}{\Pi(\lambda)}\widetilde{G}(\Omega,-\lambda,t) - \dfrac{e^{-\Omega(\lambda,t)}}{\Pi(\lambda)}\widetilde{\tilde{f}}(\Omega,-\lambda,t) & -\dfrac{i\lambda}{\Pi(\lambda)} \\ \hat{u}(\sigma_1(\lambda),t) - e^{-\Omega(\lambda,t)}\hat{u}_0(\sigma_1(\lambda)) - \dfrac{e^{-\Omega(\lambda,t)}}{\Pi(\sigma_1(\lambda))}\widetilde{G}(\Omega,\sigma_1(\lambda),t) - \dfrac{e^{-\Omega(\lambda,t)}}{\Pi(\sigma_1(\lambda))}\widetilde{\tilde{f}}(\Omega,\sigma_1(\lambda),t) & \dfrac{i\sigma_1(\lambda)}{\Pi(\sigma_1(\lambda))} \end{bmatrix}$$

and

$$(6.13) \quad e^{-\Omega(\lambda,t)}\widetilde{G}_2(\Omega(\lambda,t),t) = \frac{\Pi(\lambda)\Pi(\sigma_1(\lambda))}{i[\lambda+\sigma_1(\lambda)]} \times$$

$$\times \det\begin{bmatrix} \dfrac{1}{\Pi(\lambda)} & \hat{u}(-\lambda,t) - e^{-\Omega(\lambda,t)}\hat{u}_0(-\lambda) - \dfrac{e^{-\Omega(\lambda,t)}}{\Pi(\lambda)}\widetilde{G}(\Omega,-\lambda,t) - \dfrac{e^{-\Omega(\lambda,t)}}{\Pi(\lambda)}\widetilde{\tilde{f}}(\Omega,-\lambda,t) \\ \dfrac{1}{\Pi(\sigma_1(\lambda))} & \hat{u}(\sigma_1(\lambda),t) - e^{-\Omega(\lambda,t)}\hat{u}_0(\sigma_1(\lambda)) - \dfrac{e^{-\Omega(\lambda,t)}}{\Pi(\sigma_1(\lambda))}\widetilde{G}(\Omega,\sigma_1(\lambda),t) - \dfrac{e^{-\Omega(\lambda,t)}}{\Pi(\sigma_1(\lambda))}\widetilde{\tilde{f}}(\Omega,\sigma_1(\lambda),t) \end{bmatrix}.$$

Now, it is crucial to observe that the coefficient of $\hat{u}(-\lambda,t)$ in the quantity $e^{-\Omega(\lambda,t)}\widetilde{G}_3(\Omega(\lambda,t),t)$, as this is given by (6.12), is equal to

$$\frac{\Pi(\lambda)\Pi(\sigma_1(\lambda))}{i[\lambda+\sigma_1(\lambda)]}\frac{i\sigma_1(\lambda)}{\Pi(\sigma_1(\lambda))} = \frac{\Pi(\lambda)\sigma_1(\lambda)}{\lambda+\sigma_1(\lambda)},$$

and that

$$\int_{\mathcal{C}_1} \frac{\sigma_1(\lambda)}{\lambda + \sigma_1(\lambda)} e^{i\lambda x} \hat{u}(-\lambda, t) d\lambda = 0.$$

Similarly, we observe that the coefficient of $\hat{u}(\sigma_1(\lambda), t)$ in the quantity $e^{-\Omega(\lambda, t)} \widetilde{G}_3(\Omega(\lambda, t), t)$, as this is given by (6.12), is equal to

$$\frac{\Pi(\lambda)\Pi(\sigma_1(\lambda))}{i[\lambda + \sigma_1(\lambda)]} \frac{i\lambda}{\Pi(\lambda)} = \frac{\lambda \Pi(\sigma_1(\lambda))}{\lambda + \sigma_1(\lambda)},$$

and that

$$\int_{\mathcal{C}_1} \frac{\lambda \Pi(\sigma_1(\lambda))}{\lambda + \sigma_1(\lambda)} \frac{1}{\Pi(\lambda)} e^{i\lambda x} \hat{u}(\sigma_1(\lambda), t) d\lambda = 0,$$

where, this time, we also rely on (6.6), in order to ensure that the quotient $\Pi(\sigma_1(\lambda))/\Pi(\lambda)$ is analytic for $\lambda$ in an open neighbohood of $\lambda_1$.

Similar conclusions can be drawn also for the corresponding integrals over $\mathcal{C}_1$, which arise from the presence of these terms, namely $\hat{u}(-\lambda, t)$ and $\hat{u}(\sigma_1(\lambda), t)$, in the quantity $e^{-\Omega(\lambda, t)} \widetilde{G}_2(\Omega(\lambda, t), t)$.

Therefore, in view of (6.12) and (6.13),

(6.14) $\displaystyle\int_{\mathcal{C}_1} \frac{e^{i\lambda x - \Omega(\lambda, t)}}{\Pi(\lambda)} \widetilde{G}_3(\Omega(\lambda, t), t) d\lambda$

$$= -\int_{\mathcal{C}_1} \frac{\sigma_1(\lambda)}{\lambda + \sigma_1(\lambda)} e^{i\lambda x - \Omega(\lambda, t)} \hat{u}_0(-\lambda) d\lambda - \int_{\mathcal{C}_1} \frac{\lambda}{\lambda + \sigma_1(\lambda)} \frac{\Pi(\sigma_1(\lambda))}{\Pi(\lambda)} e^{i\lambda x - \Omega(\lambda, t)} \hat{u}_0(\sigma_1(\lambda)) d\lambda$$

$$- \int_{\mathcal{C}_1} \frac{\sigma_1(\lambda)}{\lambda + \sigma_1(\lambda)} \frac{1}{\Pi(\lambda)} e^{i\lambda x - \Omega(\lambda, t)} \widetilde{G}(\Omega(\lambda, t), -\lambda, t) d\lambda - \int_{\mathcal{C}_1} \frac{\lambda}{\lambda + \sigma_1(\lambda)} \frac{1}{\Pi(\lambda)} e^{i\lambda x - \Omega(\lambda, t)} \widetilde{G}(\Omega(\lambda, t), \sigma_1(\lambda), t) d\lambda$$

$$- \int_{\mathcal{C}_1} \frac{\sigma_1(\lambda)}{\lambda + \sigma_1(\lambda)} \frac{1}{\Pi(\lambda)} e^{i\lambda x - \Omega(\lambda, t)} \widetilde{\tilde{f}}(\Omega(\lambda, t), -\lambda, t) d\lambda - \int_{\mathcal{C}_1} \frac{\lambda}{\lambda + \sigma_1(\lambda)} \frac{1}{\Pi(\lambda)} e^{i\lambda x - \Omega(\lambda, t)} \widetilde{\tilde{f}}(\Omega(\lambda, t), \sigma_1(\lambda), t) d\lambda$$

and

(6.15) $\displaystyle\int_{\mathcal{C}_1} \frac{i\lambda e^{i\lambda x - \Omega(\lambda, t)}}{\Pi(\lambda)} \widetilde{G}_2(\Omega(\lambda, t), t) d\lambda$

$$= \int_{\mathcal{C}_1} \frac{\lambda}{\lambda + \sigma_1(\lambda)} e^{i\lambda x - \Omega(\lambda, t)} \hat{u}_0(-\lambda) d\lambda - \int_{\mathcal{C}_1} \frac{\lambda}{\lambda + \sigma_1(\lambda)} \frac{\Pi(\sigma_1(\lambda))}{\Pi(\lambda)} e^{i\lambda x - \Omega(\lambda, t)} \hat{u}_0(\sigma_1(\lambda)) d\lambda$$

$$+ \int_{\mathcal{C}_1} \frac{\lambda}{\lambda + \sigma_1(\lambda)} \frac{1}{\Pi(\lambda)} e^{i\lambda x - \Omega(\lambda, t)} \widetilde{G}(\Omega(\lambda, t), -\lambda, t) d\lambda - \int_{\mathcal{C}_1} \frac{\lambda}{\lambda + \sigma_1(\lambda)} \frac{1}{\Pi(\lambda)} e^{i\lambda x - \Omega(\lambda, t)} \widetilde{G}(\Omega(\lambda, t), \sigma_1(\lambda), t) d\lambda$$

$$+ \int_{\mathcal{C}_1} \frac{\lambda}{\lambda + \sigma_1(\lambda)} \frac{1}{\Pi(\lambda)} e^{i\lambda x - \Omega(\lambda, t)} \widetilde{\tilde{f}}(\Omega(\lambda, t), -\lambda, t) d\lambda - \int_{\mathcal{C}_1} \frac{\lambda}{\lambda + \sigma_1(\lambda)} \frac{1}{\Pi(\lambda)} e^{i\lambda x - \Omega(\lambda, t)} \widetilde{\tilde{f}}(\Omega(\lambda, t), \sigma_1(\lambda), t) d\lambda.$$

Similarly, solving the system of equations (6.8) & (6.10), for $\lambda \in \mathcal{C}_2$, we compute

(6.16) $\displaystyle\int_{\mathcal{C}_2} \frac{e^{i\lambda x - \Omega(\lambda, t)}}{\Pi(\lambda)} \widetilde{G}_3(\Omega(\lambda, t), t) d\lambda$

$\qquad\qquad$ = the *RHS* of (6.14), where, *however*, $\mathcal{C}_1$ is replaced by $\mathcal{C}_2$ and $\sigma_1(\lambda)$ is replaced by $\sigma_2(\lambda)$,

and

(6.17) $\displaystyle\int_{\mathcal{C}_2} \frac{i\lambda e^{i\lambda x - \Omega(\lambda, t)}}{\Pi(\lambda)} \widetilde{G}_2(\Omega(\lambda, t), t) d\lambda$

= the *RHS* of (6.15), where, *however*, $C_1$ is replaced by $C_2$ and $\sigma_1(\lambda)$ by is replaced by $\sigma_2(\lambda)$.

Finally, substituting the integrals (6.14), (6.15), (6.16) and (6.17) in (6.7), we obtain the integral representation of the solution of problem (6.1), in the case $4\alpha > \beta^2$.

**The case $4\alpha < \beta^2$ with $\beta > 0$.** In this case the zeros of the polynomial $\Pi(\lambda) = 1 + \beta\lambda^2 + \alpha\lambda^4$ are the numbers $\lambda_1$, $\lambda_2$, $-\lambda_1$, $-\lambda_2$, where

$$\lambda_1 = i\sqrt{\frac{\beta + \sqrt{\beta^2 - 4\alpha}}{2\alpha}} \quad \text{and} \quad \lambda_2 = i\sqrt{\frac{\beta - \sqrt{\beta^2 - 4\alpha}}{2\alpha}},$$

and the construction of the solution of problem (6.1) is similar to the previous case.

**The case $4\alpha = \beta^2$** It seems that in this case we do not have sufficient data for problem (6.1), since the polynomial $\Pi(\lambda) = 1 + \beta\lambda^2 + \alpha\lambda^4$ has two double zeros. For example, if we include in the data of the problem the datum $u_{xx}(0,t) = g_2(t)$, then it seems that we can carry out the construction of the solution.

## 7. The operator $\partial_t + \alpha \partial_{xxxxt} - \beta \partial_{xxt} - \gamma(t)\partial_x$

**Problem 6** Solve

(7.1)
$$\begin{cases} \partial_t u = -\alpha \partial_{xxxxt} u + \beta \partial_{xxt} u - \gamma(t) \partial_x u + f(x,t), \ (x,t) \in \mathbb{R}^+ \times \mathbb{R}^+, \\ u(x,0) = u_0(x), \ x \in \mathbb{R}^+, \\ u(0,t) = g_0(t), \ t \in \mathbb{R}^+, \\ u_x(0,t) = g_1(t), \ t \in \mathbb{R}^+, \end{cases}$$

for $u(x,t)$.

**Derivation of the solution** For $\lambda \in \mathbb{C}$ with $\text{Im}\,\lambda \leq 0$, the differential equation in (7.1) leads to

$$(1 + \beta\lambda^2 + \alpha\lambda^4)\frac{\partial}{\partial t}[\hat{u}(\lambda,t)] + i\lambda\gamma(t)\hat{u}(\lambda,t) = \alpha[g'_3(t) + i\lambda g'_2(t) - \lambda^2 g'_1(t) - i\lambda^3 g'_0(t)]$$

$$- \beta[g'_1(t) + i\lambda g'_0(t)] + i\lambda\gamma(t)g_0(t) + \hat{f}(\lambda,t)$$

$\Rightarrow$

$$\frac{\partial}{\partial t}[\hat{u}(\lambda,t)] + \frac{\gamma(t)i\lambda}{1 + \beta\lambda^2 + \alpha\lambda^4}\hat{u}(\lambda,t) = \frac{\alpha}{1 + \beta\lambda^2 + \alpha\lambda^4}[g'_3(t) + i\lambda g'_2(t) - \lambda^2 g'_1(t) - i\lambda^3 g'_0(t)]$$

$$- \frac{\beta}{1 + \beta\lambda^2 + \alpha\lambda^4}[g'_1(t) + i\lambda g'_0(t)] + \frac{\gamma(t)i\lambda}{1 + \beta\lambda^2 + \alpha\lambda^4}g_0(t) + \frac{\hat{f}(\lambda,t)}{1 + \beta\lambda^2 + \alpha\lambda^4}$$

$\Rightarrow$

(7.2) $\frac{\partial}{\partial t}[\hat{u}(\lambda,t)] + \frac{\gamma(t)i\lambda}{1 + \beta\lambda^2 + \alpha\lambda^4}\hat{u}(\lambda,t)$

$$= \frac{\alpha}{1 + \beta\lambda^2 + \alpha\lambda^4}g'_3(t) + \frac{\alpha i\lambda}{1 + \beta\lambda^2 + \alpha\lambda^4}g'_2(t) + \frac{G(\lambda,t)}{1 + \beta\lambda^2 + \alpha\lambda^4} + \frac{\hat{f}(\lambda,t)}{1 + \beta\lambda^2 + \alpha\lambda^4},$$

where we have set

$$G(\lambda,t) := -\alpha[\lambda^2 g'_1(t) + i\lambda^3 g'_0(t)] - \beta[g'_1(t) + i\lambda g'_0(t)] + \gamma(t)i\lambda g_0(t)$$

$$= -(\alpha\lambda^2 + \beta)g'_1(t) - i\lambda(\alpha\lambda^2 + \beta)g'_0(t) + \gamma(t)i\lambda g_0(t).$$

Furthermore, setting

$$\omega(\lambda,t) = \frac{\gamma(t)i\lambda}{1+\beta\lambda^2+\alpha\lambda^4}, \quad \Omega(\lambda,t) = \int_{\tau=0}^{t}\omega(\lambda,\tau)d\tau = \int_{\tau=0}^{t}\frac{\gamma(\tau)i\lambda}{1+\beta\lambda^2+\alpha\lambda^4}d\tau = \frac{i\lambda}{1+\beta\lambda^2+\alpha\lambda^4}\int_{\tau=0}^{t}\gamma(\tau)d\tau,$$

we write (7.2) in the following way:

$$\frac{\partial}{\partial t}[e^{\Omega(\lambda,t)}\hat{u}(\lambda,t)] = \frac{\alpha e^{\Omega(\lambda,t)}g'_3(t)}{1+\beta\lambda^2+\alpha\lambda^4} + \frac{\alpha i\lambda e^{\Omega(\lambda,t)}g'_2(t)}{1+\beta\lambda^2+\alpha\lambda^4} - \frac{e^{\Omega(\lambda,t)}G(\lambda,t)}{1+\beta\lambda^2+\alpha\lambda^4} + \frac{e^{\Omega(\lambda,t)}\hat{f}(\lambda,t)}{1+\beta\lambda^2+\alpha\lambda^4}.$$

Integrating the above equation, we find

$$(7.3) \quad e^{\Omega(\lambda,t)}\hat{u}(\lambda,t) = \hat{u}_0(\lambda) + \alpha\int_{\tau=0}^{t}\frac{e^{\Omega(\lambda,\tau)}g'_3(\tau)}{1+\beta\lambda^2+\alpha\lambda^4}d\tau + \alpha\int_{\tau=0}^{t}\frac{i\lambda e^{\Omega(\lambda,\tau)}g'_2(\tau)}{1+\beta\lambda^2+\alpha\lambda^4}d\tau$$

$$+ \int_{\tau=0}^{t}\frac{e^{\Omega(\lambda,\tau)}G(\lambda,\tau)}{1+\beta\lambda^2+\alpha\lambda^4}d\tau + \int_{\tau=0}^{t}\frac{e^{\Omega(\lambda,\tau)}\hat{f}(\lambda,\tau)}{1+\beta\lambda^2+\alpha\lambda^4}d\tau.$$

Defining

$$(g'_3\,\widetilde{\,})(\Omega(\lambda,t),t) = \int_{\tau=0}^{t}e^{\Omega(\lambda,\tau)}g'_3(\tau)d\tau, \quad (g'_2\,\widetilde{\,})(\Omega(\lambda,t),t) = \int_{\tau=0}^{t}e^{\Omega(\lambda,\tau)}g'_2(\tau)d\tau,$$

$$\widetilde{G}(\Omega(\lambda,t),\lambda,t) = \int_{\tau=0}^{t}e^{\Omega(\lambda,\tau)}G(\lambda,\tau)d\tau \text{ and } \widetilde{\hat{f}}(\Omega(\lambda,t),\lambda,t) = \int_{\tau=0}^{t}e^{\Omega(\lambda,\tau)}\hat{f}(\lambda,\tau)d\tau,$$

we write (7.3) as follows:

$$(7.4) \quad \hat{u}(\lambda,t) = e^{-\Omega(\lambda,t)}\hat{u}_0(\lambda) + \alpha\frac{e^{-\Omega(\lambda,t)}}{1+\beta\lambda^2+\alpha\lambda^4}(g'_3\,\widetilde{\,})(\Omega(\lambda,t),t) + \alpha\frac{i\lambda e^{-\Omega(\lambda,t)}}{1+\beta\lambda^2+\alpha\lambda^4}(g'_2\,\widetilde{\,})(\Omega(\lambda,t),t)$$

$$+ \frac{e^{-\Omega(\lambda,t)}}{1+\beta\lambda^2+\alpha\lambda^4}\widetilde{G}(\Omega(\lambda,t),\lambda,t) + \frac{e^{-\Omega(\lambda,t)}}{1+\beta\lambda^2+\alpha\lambda^4}\widetilde{\hat{f}}(\Omega(\lambda,t),\lambda,t), \text{ for } \lambda\in\mathbb{C} \text{ with } \operatorname{Im}\lambda\leq 0.$$

Now (7.4) gives

$$(7.5) \quad 2\pi u(x,t) = \int_{-\infty}^{\infty}e^{i\lambda x}\hat{u}(\lambda,t)d\lambda$$

$$= \int_{-\infty}^{\infty}e^{i\lambda x-\Omega(\lambda,t)}\hat{u}_0(\lambda)d\lambda + \alpha\int_{-\infty}^{\infty}\frac{e^{i\lambda x-\Omega(\lambda,t)}}{1+\beta\lambda^2+\alpha\lambda^4}(g'_3\,\widetilde{\,})(\Omega(\lambda,t),t)d\lambda + \alpha\int_{-\infty}^{\infty}\frac{i\lambda e^{i\lambda x-\Omega(\lambda,t)}}{1+\beta\lambda^2+\alpha\lambda^4}(g'_2\,\widetilde{\,})(\Omega(\lambda,t),t)d\lambda$$

$$+ \int_{-\infty}^{\infty}\frac{e^{i\lambda x-\Omega(\lambda,t)}}{1+\beta\lambda^2+\alpha\lambda^4}\widetilde{G}(\Omega(\lambda,t),\lambda,t)d\lambda + \int_{-\infty}^{\infty}\frac{e^{i\lambda x-\Omega(\lambda,t)}}{1+\beta\lambda^2+\alpha\lambda^4}\widetilde{\hat{f}}(\Omega(\lambda,t),\lambda,t)d\lambda.$$

We continue the construction, assuming $4\alpha > \beta^2$.

**The case $4\alpha > \beta^2$** We consider the polynomial $\Pi(\lambda) := 1+\beta\lambda^2+\alpha\lambda^4$ and its zeros $\lambda_1, \lambda_2, -\lambda_1, -\lambda_2$, as in section 6. Considering the contours $C_1$ and $C_2$, which were defined in section 6, we write (7.5) as follows:

$$(7.6) \quad 2\pi u(x,t) = \int_{-\infty}^{\infty}e^{i\lambda x-\Omega(\lambda,t)}\hat{u}_0(\lambda)d\lambda + \alpha\int_{C_1}\frac{e^{i\lambda x-\Omega(\lambda,t)}}{\Pi(\lambda)}(g'_3\,\widetilde{\,})(\Omega(\lambda,t),t)d\lambda + \alpha\int_{C_2}\frac{e^{i\lambda x-\Omega(\lambda,t)}}{\Pi(\lambda)}(g'_3\,\widetilde{\,})(\Omega(\lambda,t),t)d\lambda$$

$$+ \alpha\int_{C_1}\frac{i\lambda e^{i\lambda x-\Omega(\lambda,t)}}{\Pi(\lambda)}(g'_2\,\widetilde{\,})(\Omega(\lambda,t),t)d\lambda + \alpha\int_{C_2}\frac{i\lambda e^{i\lambda x-\Omega(\lambda,t)}}{\Pi(\lambda)}(g'_2\,\widetilde{\,})(\Omega(\lambda,t),t)d\lambda$$

$$+ \int_{C_1} \frac{e^{i\lambda x - \Omega(\lambda,t)}}{\Pi(\lambda)} \widetilde{G}(\Omega(\lambda,t),\lambda,t) d\lambda + \int_{C_2} \frac{e^{i\lambda x - \Omega(\lambda,t)}}{\Pi(\lambda)} \widetilde{G}(\Omega(\lambda,t),\lambda,t) d\lambda + \int_{-\infty}^{\infty} \frac{e^{i\lambda x - \Omega(\lambda,t)}}{\Pi(\lambda)} \widetilde{f}(\Omega(\lambda,t),\lambda,t) d\lambda.$$

Next we consider the equation $\omega(\sigma(\lambda),t) = \omega(\lambda,t)$, i.e.,

$$\frac{\sigma(\lambda)}{1+\beta[\sigma(\lambda)]^2 + \alpha[\sigma(\lambda)]^4} = \frac{\lambda}{1+\beta\lambda^2 + \alpha\lambda^4} \Rightarrow \lambda\{1+\beta[\sigma(\lambda)]^2 + \alpha[\sigma(\lambda)]^4\} = (1+\beta\lambda^2 + \alpha\lambda^4)\sigma(\lambda).$$

This leads to the algebraic equation $\alpha\lambda[\sigma(\lambda)]^4 + \beta\lambda[\sigma(\lambda)]^2 - (1+\beta\lambda^2 + \alpha\lambda^4)\sigma(\lambda) + \lambda = 0$, for the function $\sigma(\lambda)$. We consider the polynomial

$$Q(\lambda,\sigma) := \alpha\lambda\sigma^4 + \beta\lambda\sigma^2 - (1+\beta\lambda^2 + \alpha\lambda^4)\sigma + \lambda$$
$$= (\sigma - \lambda)[\alpha\lambda\sigma^3 + \alpha\lambda^2\sigma^2 + (\beta\lambda + \alpha\lambda^3)\sigma - 1],$$

in the two complex variables, $\lambda$ and $\sigma$, i.e., $Q(\lambda,\sigma) \in \mathbb{C}[\lambda,\sigma]$.

Fixing a $j \in \{1,2\}$, let us observe that the zeros of the polynomial $Q(\lambda_j,\sigma) = Q(\lambda,\sigma)|_{\lambda=\lambda_j} \in \mathbb{C}[\sigma]$, considered as a polynomial of $\sigma$, are the numbers $\lambda_1, \lambda_2, -\lambda_1, -\lambda_2$. In particular, these zeros are simple and, therefore,

(7.7) $$Q_\sigma(\lambda_j,\sigma) = \frac{\partial Q(\lambda,\sigma)}{\partial \sigma}\bigg|_{\lambda=\lambda_j} \neq 0, \text{ for } \sigma \in \{\lambda_1, \lambda_2, -\lambda_1, -\lambda_2\}.$$

Applying (7.7) with $j=1$ and $\sigma=-\lambda_1$, we obtain

$$\frac{\partial Q(\lambda,\sigma)}{\partial \sigma}\bigg|_{(\lambda,\sigma)=(\lambda_1,-\lambda_1)} \neq 0.$$

Also $Q(\lambda,\sigma)|_{(\lambda,\sigma)=(\lambda_1,-\lambda_1)} = 0$, and, therefore, by the implicit function theorem, there exists a function $\sigma_{1,1}(\lambda)$, analytic for $\lambda$ in an open neighborhood of $\lambda_1$, such that

(7.8) $$\sigma_{1,1}(\lambda_1) = -\lambda_1 \text{ and } Q(\lambda,\sigma_{1,1}(\lambda)) = 0, \text{ for } \lambda \text{ close to } \lambda_1.$$

Moreover,

$$Q_\lambda(\lambda,\sigma_{1,1}(\lambda)) + Q_\sigma(\lambda,\sigma_{1,1}(\lambda))\frac{d\sigma_{1,1}(\lambda)}{d\lambda} = 0 \Rightarrow \frac{d\sigma_{1,1}(\lambda)}{d\lambda}\bigg|_{\lambda=\lambda_1} = -\frac{Q_\lambda(\lambda,\sigma_{1,1}(\lambda))}{Q_\sigma(\lambda,\sigma_{1,1}(\lambda))}\bigg|_{\lambda=\lambda_1} = -\frac{Q_\lambda(\lambda_1,-\lambda_1)}{Q_\sigma(\lambda_1,-\lambda_1)} \neq 0,$$

and, therefore, the function $\sigma_{1,1}(\lambda)$ is biholomorphic from an open neighborhood of $\lambda_1$ to an open neighborhood of $-\lambda_1$. Also $\sigma_{1,1}(\lambda)$ satisfies

(7.9) $$\omega(\sigma_{1,1}(\lambda),t) = \omega(\lambda,t) \text{ and } \Omega(\sigma_{1,1}(\lambda),t) = \Omega(\lambda,t), \text{ for } \lambda \text{ close to } \lambda_1.$$

Similarly, $Q(\lambda,\sigma)|_{(\lambda,\sigma)=(\lambda_1,-\lambda_2)} = 0$ and, as it follows from (7.7),

$$\frac{\partial Q(\lambda,\sigma)}{\partial \sigma}\bigg|_{(\lambda,\sigma)=(\lambda_1,-\lambda_2)} \neq 0.$$

Therefore, there exists a function $\sigma_{1,2}(\lambda)$, analytic for $\lambda$ in an open neighborhood of $\lambda_1$, such that

$$\sigma_{1,2}(\lambda_1) = -\lambda_2 \text{ and } Q(\lambda,\sigma_{1,2}(\lambda)) = 0, \text{ for } \lambda \text{ close to } \lambda_1.$$

Also, the function $\sigma_{1,2}(\lambda)$ is biholomorphic from an open neighborhood of $\lambda_1$ to an open neighborhood of $-\lambda_2$ and satisfies

(7.10) $\quad \omega(\sigma_{1,2}(\lambda),t) = \omega(\lambda,t)$ and $\Omega(\sigma_{1,2}(\lambda),t) = \Omega(\lambda,t)$, for $\lambda$ close to $\lambda_1$.

Working with $\lambda$ close to $\lambda_2$, we can also construct analytic functions $\sigma_{2,1}(\lambda)$ and $\sigma_{2,2}(\lambda)$, defined both on an open neighborhood of $\lambda_2$, such that

$$\sigma_{2,1}(\lambda_2) = -\lambda_1, \; \sigma_{2,2}(\lambda_2) = -\lambda_2, \; \left.\frac{d\sigma_{2,1}(\lambda)}{d\lambda}\right|_{\lambda=\lambda_2} \neq 0, \; \left.\frac{d\sigma_{2,2}(\lambda)}{d\lambda}\right|_{\lambda=\lambda_2} \neq 0,$$

and, for $\lambda$ close to $\lambda_2$,

(7.11) $\quad \omega(\sigma_{2,1}(\lambda),t) = \omega(\lambda,t)$, $\Omega(\sigma_{2,1}(\lambda),t) = \Omega(\lambda,t)$, $\omega(\sigma_{2,2}(\lambda),t) = \omega(\lambda,t)$, $\Omega(\sigma_{2,2}(\lambda),t) = \Omega(\lambda,t)$.

*Note* The functions $\sigma_{1,1}(\lambda)$ and $\sigma_{1,2}(\lambda)$, are solutions of the cubic equation

(7.12) $\quad \alpha\lambda[\sigma(\lambda)]^3 + \alpha\lambda^2[\sigma(\lambda)]^2 + (\beta\lambda + \alpha\lambda^3)\sigma(\lambda) - 1 = 0$,

defined for $\lambda$ close to $\lambda_1$, with the specified properties, while the functions $\sigma_{2,1}(\lambda)$ and $\sigma_{2,2}(\lambda)$, are solutions of (7.12), defined for $\lambda$ close to $\lambda_2$. In particular, these functions can be explicitly expressed in terms of the parameters $\alpha$ and $\beta$, using cubic and square roots, via Cardano's formula. Our efforts, in the previous paragraph, were aiming to show that appropriate choices of the cubic and square roots in Cardano's formula, give solutions of (7.12) with properties, which allow us to continue the construction of the solution of problem (7.1).

In view of the properties of the functions $\sigma_{1,1}(\lambda)$, $\sigma_{1,2}(\lambda)$, $\sigma_{2,1}(\lambda)$, $\sigma_{2,2}(\lambda)$, we see that (7.4) implies:

(7.13) $\quad \hat{u}(\sigma_{1,j}(\lambda),t) = e^{-\Omega(\lambda,t)}\hat{u}_0(\sigma_{1,j}(\lambda)) + \alpha\frac{e^{-\Omega(\lambda,t)}}{\Pi(\sigma_{1,j}(\lambda))}(g'_3\widetilde{\;})(\Omega(\lambda,t),t) + \alpha\frac{i\sigma_{1,j}(\lambda)e^{-\Omega(\lambda,t)}}{\Pi(\sigma_{1,j}(\lambda))}(g'_2\widetilde{\;})(\Omega(\lambda,t),t)$

$\quad + \frac{e^{-\Omega(\lambda,t)}}{\Pi(\sigma_{1,j}(\lambda))}\widetilde{G}(\Omega(\lambda,t),\sigma_{1,j}(\lambda),t) + \frac{e^{-\Omega(\lambda,t)}}{\Pi(\sigma_{1,j}(\lambda))}\widetilde{\tilde{f}}(\Omega(\lambda,t),\sigma_{1,j}(\lambda),t)$, for $\lambda \in C_1$ and $j \in \{1,2\}$.

and

(7.14) $\quad \hat{u}(\sigma_{2,j}(\lambda),t) = e^{-\Omega(\lambda,t)}\hat{u}_0(\sigma_{2,j}(\lambda)) + \alpha\frac{e^{-\Omega(\lambda,t)}}{\Pi(\sigma_{2,j}(\lambda))}(g'_3\widetilde{\;})(\Omega(\lambda,t),t) + \alpha\frac{i\sigma_{2,j}(\lambda)e^{-\Omega(\lambda,t)}}{\Pi(\sigma_{2,j}(\lambda))}(g'_2\widetilde{\;})(\Omega(\lambda,t),t)$

$\quad + \frac{e^{-\Omega(\lambda,t)}}{\Pi(\sigma_{2,j}(\lambda))}\widetilde{G}(\Omega(\lambda,t),\sigma_{2,j}(\lambda),t) + \frac{e^{-\Omega(\lambda,t)}}{\Pi(\sigma_{2,j}(\lambda))}\widetilde{\tilde{f}}(\Omega(\lambda,t),\sigma_{2,j}(\lambda),t)$, for $\lambda \in C_2$ and $j \in \{1,2\}$.

We consider (7.13) as a system of algebraic equations in the unknown quantities

(7.15) $\quad e^{-\Omega(\lambda,t)}(g'_3\widetilde{\;})(\Omega(\lambda,t),t)$ and $e^{-\Omega(\lambda,t)}(g'_2\widetilde{\;})(\Omega(\lambda,t),t)$, for $\lambda \in C_1$.

The determinant of the unknowns (7.15), in this system, is

$$\det\begin{bmatrix} \dfrac{\alpha}{\Pi(\sigma_{1,1}(\lambda))} & \dfrac{\alpha i\sigma_{1,1}(\lambda)}{\Pi(\sigma_{1,1}(\lambda))} \\ \dfrac{\alpha}{\Pi(\sigma_{1,2}(\lambda))} & \dfrac{\alpha i\sigma_{1,2}(\lambda)}{\Pi(\sigma_{1,2}(\lambda))} \end{bmatrix} = \frac{\alpha^2 i[\sigma_{1,2}(\lambda) - \sigma_{1,1}(\lambda)]}{\Pi(\sigma_{1,1}(\lambda))\Pi(\sigma_{1,2}(\lambda))} \neq 0, \text{ for } \lambda \text{ close to } \lambda_1.$$

An analogous relation holds also for the system (7.14), and the unknowns (7.15), this time, however, for $\lambda$ close to $\lambda_2$.

Another crucial point, that we have to check, is the vanishing of all the integrals which will appear in (7.6) and which will contain the quantities $\hat{u}(\sigma_{1,1}(\lambda),t)$, $\hat{u}(\sigma_{1,2}(\lambda),t)$, $\hat{u}(\sigma_{2,1}(\lambda),t)$, $\hat{u}(\sigma_{2,1}(\lambda),t)$, when we will substitute the values of (7.15) (for both $\lambda \in \mathcal{C}_1$ and $\lambda \in \mathcal{C}_2$) in (7.6), as these will be computed by solving the systems (7.14) and (7.15). For example, we have that

$$(7.16) \quad \int_{\mathcal{C}_1} e^{i\lambda x} \hat{u}(\sigma_{1,1}(\lambda),t) \frac{\alpha i \sigma_{1,2}(\lambda)}{\Pi(\sigma_{1,2}(\lambda))} \left\{ \frac{\alpha^2 i [\sigma_{1,2}(\lambda) - \sigma_{1,1}(\lambda)]}{\Pi(\sigma_{1,1}(\lambda))\Pi(\sigma_{1,2}(\lambda))} \right\}^{-1} \frac{1}{\Pi(\lambda)} d\lambda = 0,$$

and this follows from the analyticity of the integrand, for $\lambda$ close to $\lambda_1$, whose analyticity, in turn, follows from the equation

$$\frac{\Pi(\sigma_{1,1}(\lambda))}{\Pi(\lambda)} = \frac{\sigma_{1,1}(\lambda)}{\lambda}.$$

Thus, solving the systems (7.14) and (7.15) and substituting the valus of (7.15) (for both $\lambda \in \mathcal{C}_1$ and $\lambda \in \mathcal{C}_2$), in (7.6), we obtain the formula for the solution of problem (7.1), in the case $4\alpha > \beta^2$.

*Note* The comments that we made at the end of section 6, for problem 5, concerning the cases $4\alpha < \beta^2$ and $4\alpha = \beta^2$, apply also in the case of problem 6.

## 8. The operator $\partial_t - \alpha \partial_{xxxxxxt} + \beta \partial_{xxxxt} - \gamma \partial_{xxt} - \delta(t) \partial_{xx}$

**Problem 7** Solve

$$(8.1) \quad \begin{cases} \partial_t u = \alpha \partial_{xxxxxxt} u - \beta \partial_{xxxxt} u + \gamma \partial_{xxt} u + \delta(t) \partial_{xx} u + f(x,t), \ (x,t) \in \mathbb{R}^+ \times \mathbb{R}^+, \\ u(x,0) = u_0(x), \ x \in \mathbb{R}^+, \\ u(0,t) = g_0(t), \ t \in \mathbb{R}^+, \\ u_x(0,t) = g_1(t), \ t \in \mathbb{R}^+, \\ u_{xx}(0,t) = g_2(t), \ t \in \mathbb{R}^+, \end{cases}$$

for $u(x,t)$.

***Derivation of the solution*** For $\lambda \in \mathbb{C}$ with $\operatorname{Im}\lambda \leq 0$, the differential equation in (8.1) leads to

$$(1+\gamma\lambda^2 + \beta\lambda^4 + \alpha\lambda^6) \frac{\partial}{\partial t}[\hat{u}(\lambda,t)] + \delta(t)\lambda^2 \hat{u}(\lambda,t)$$
$$= -\alpha[g'_5(t) + i\lambda g'_4(t) + (i\lambda)^2 g'_3(t) + (i\lambda)^3 g'_2(t) + (i\lambda)^4 g'_1(t) + (i\lambda)^5 g'_0(t)]$$
$$+ \beta[g'_3(t) + i\lambda g'_2(t) + (i\lambda)^2 g'_1(t) + (i\lambda)^3 g'_0(t)]$$
$$- \gamma[g'_1(t) + i\lambda g'_0(t)] - \delta(t)[g_1(t) + i\lambda g_0(t)] + \hat{f}(\lambda,t)$$

$\Rightarrow$

$$(8.2) \quad \frac{\partial}{\partial t}[\hat{u}(\lambda,t)] + \frac{\delta(t)\lambda^2}{\Pi(\lambda)} \hat{u}(\lambda,t)$$
$$= -\frac{\alpha}{\Pi(\lambda)}[g'_5(t) + i\lambda g'_4(t) + (i\lambda)^2 g'_3(t) + (i\lambda)^3 g'_2(t) + (i\lambda)^4 g'_1(t) + (i\lambda)^5 g'_0(t)]$$
$$+ \frac{\beta}{\Pi(\lambda)}[g'_3(t) + i\lambda g'_2(t) + (i\lambda)^2 g'_1(t) + (i\lambda)^3 g'_0(t)]$$
$$- \frac{\gamma}{\Pi(\lambda)}[g'_1(t) + i\lambda g'_0(t)] - \frac{\delta(t)}{\Pi(\lambda)}[g_1(t) + i\lambda g_0(t)] + \frac{\hat{f}(\lambda,t)}{\Pi(\lambda)},$$

where we have set $\Pi(\lambda) := 1 + \gamma\lambda^2 + \beta\lambda^4 + \alpha\lambda^6$.

Next, setting

$$G(\lambda,t) = -\alpha[(i\lambda)^3 g'_2(t) + (i\lambda)^4 g'_1(t) + (i\lambda)^5 g'_0(t)]$$

$$+ i\lambda g'_2(t) + (i\lambda)^2 g'_1(t) + (i\lambda)^3 g'_0(t)] - \gamma[g'_1(t) + i\lambda g'_0(t)] - \delta(t)[g_1(t) + i\lambda g_0(t)]$$

$$= i\lambda(\alpha\lambda^2 + \beta)g'_2(t) - (\alpha\lambda^4 + \beta\lambda^2 + \gamma)g'_1(t) - i\lambda(\alpha\lambda^4 + \beta\lambda^2 + \gamma)g'_0(t) - \delta(t)[g_1(t) + i\lambda g_0(t)],$$

we write (8.2) as follows:

$$(8.3) \quad \frac{\partial}{\partial t}[\hat{u}(\lambda,t)] + \frac{\delta(t)\lambda^2}{\Pi(\lambda)}\hat{u}(\lambda,t) = -\frac{\alpha g'_5(t)}{\Pi(\lambda)} - \frac{\alpha i\lambda g'_4(t)}{\Pi(\lambda)} + \frac{(\alpha\lambda^2 + \beta)g'_3(t)}{\Pi(\lambda)} + \frac{G(\lambda,t)}{\Pi(\lambda)} + \frac{\hat{f}(\lambda,t)}{\Pi(\lambda)}.$$

Furthermore, setting

$$\omega(\lambda,t) = \frac{\delta(t)\lambda^2}{1 + \gamma\lambda^2 + \beta\lambda^4 + \alpha\lambda^6} = \frac{\delta(t)\lambda^2}{\Pi(\lambda)} \quad \text{and} \quad \Omega(\lambda,t) = \int_{\tau=0}^{t} \omega(\lambda,\tau)d\tau = \frac{\lambda^2}{\Pi(\lambda)}\int_{\tau=0}^{t}\delta(\tau)d\tau,$$

we write (8.3) in the following way:

$$\frac{\partial}{\partial t}[e^{\Omega(\lambda,t)}\hat{u}(\lambda,t)] = -\frac{\alpha e^{\Omega(\lambda,t)}g'_5(t)}{\Pi(\lambda)} - \frac{\alpha i\lambda e^{\Omega(\lambda,t)}g'_4(t)}{\Pi(\lambda)} + \frac{(\alpha\lambda^2 + \beta)e^{\Omega(\lambda,t)}g'_3(t)}{\Pi(\lambda)} + \frac{e^{\Omega(\lambda,t)}G(\lambda,t)}{\Pi(\lambda)} + \frac{e^{\Omega(\lambda,t)}\hat{f}(\lambda,t)}{\Pi(\lambda)}.$$

Integrating, we obtain

$$(8.4) \quad \hat{u}(\lambda,t) = \hat{u}_0(\lambda) - \frac{\alpha e^{-\Omega(\lambda,t)}(\widetilde{g'_5})(\Omega(\lambda,t),t)}{\Pi(\lambda)} - \frac{\alpha i\lambda e^{-\Omega(\lambda,t)}(\widetilde{g'_4})(\Omega(\lambda,t),t)}{\Pi(\lambda)}$$

$$+ \frac{(\alpha\lambda^2 + \beta)e^{-\Omega(\lambda,t)}(\widetilde{g'_3})(\Omega(\lambda,t),t)}{\Pi(\lambda)} + \frac{e^{-\Omega(\lambda,t)}\widetilde{G}(\Omega(\lambda,t),\lambda,t)}{\Pi(\lambda)} + \frac{e^{-\Omega(\lambda,t)}\widetilde{\hat{f}}(\Omega(\lambda,t),\lambda,t)}{\Pi(\lambda)},$$

where

$$(\widetilde{g'_j})(\Omega(\lambda,t),t) := \int_{\tau=0}^{t} e^{\Omega(\lambda,\tau)} g'_j(\tau)d\tau, \; j \in \{5, 4, 3\},$$

$$\widetilde{G}(\Omega(\lambda,t),\lambda,t) := \int_{\tau=0}^{t} e^{\Omega(\lambda,\tau)} G(\lambda,\tau)d\tau \quad \text{and} \quad \widetilde{\hat{f}}(\Omega(\lambda,t),\lambda,t) := \int_{\tau=0}^{t} e^{\Omega(\lambda,\tau)} \hat{f}(\lambda,\tau)d\tau.$$

*Assumpion on the parameters $\alpha, \beta, \gamma$.* We assume that the polynomial $\varpi(\mu) := 1 + \gamma\mu + \beta\mu^2 + \alpha\mu^3$ has dinstict zeros, i.e., the resultant $\mathcal{R}(\varpi, \varpi')$, of the polynomials $\varpi(\mu)$ and $\varpi'(\mu) = \gamma + 2\beta\mu + 3\alpha\mu^2$, is not zero. More precisely

$$(8.5) \quad \mathcal{R}(\varpi,\varpi') = \det\begin{bmatrix} \alpha & 0 & 3\alpha & 0 & 0 \\ \beta & \alpha & 2\beta & 3\alpha & 0 \\ \gamma & \beta & \gamma & 2\beta & 3\alpha \\ 1 & \gamma & 0 & \gamma & 2\beta \\ 0 & 1 & 0 & 0 & \gamma \end{bmatrix} \neq 0.$$

It follows from (8.5) that the zeros of the polynomial $\Pi(\lambda) = 1 + \gamma\lambda^2 + \beta\lambda^4 + \alpha\lambda^6$ are distinct, and in fact three of them lie in the upper half-plane $\{\text{Im}\,\lambda > 0\}$ and the other three lie in the lower half-plane $\{\text{Im}\,\lambda < 0\}$. Let $\lambda_1, \lambda_2, \lambda_3$ be the zeros of $\Pi(\lambda)$ with $\text{Im}\,\lambda_1 > 0, \text{Im}\,\lambda_2 > 0, \text{Im}\,\lambda_3 > 0$, in which case the zeros of $\Pi(\lambda)$ in the half-plane $\{\text{Im}\,\lambda < 0\}$ are the numbers $-\lambda_1, -\lambda_2, -\lambda_3$.

Next, we consider the contours $\mathcal{C}(\lambda_1), \mathcal{C}(\lambda_2), \mathcal{C}(\lambda_3)$, to be simple closed curves, around the points $\lambda_1, \lambda_2, \lambda_3$, sufficiently small and oriented counterclockwise. Then (8.4) gives the following equation

$$(8.6)\quad 2\pi u(x,t) = \int_{-\infty}^{\infty} e^{i\lambda x}\hat{u}(\lambda,t)d\lambda = \int_{-\infty}^{\infty} e^{i\lambda x}\hat{u}_0(\lambda)d\lambda - \int_{\mathcal{C}(\lambda_1)+\mathcal{C}(\lambda_2)+\mathcal{C}(\lambda_3)} \frac{\alpha e^{i\lambda x-\Omega(\lambda,t)}(g'_5\,\widetilde{)}(\Omega(\lambda,t),t)}{\Pi(\lambda)}d\lambda$$

$$-\int_{\mathcal{C}(\lambda_1)+\mathcal{C}(\lambda_2)+\mathcal{C}(\lambda_3)} \frac{\alpha i\lambda e^{i\lambda x-\Omega(\lambda,t)}(g'_4\,\widetilde{)}(\Omega(\lambda,t),t)}{\Pi(\lambda)}d\lambda + \int_{\mathcal{C}(\lambda_1)+\mathcal{C}(\lambda_2)+\mathcal{C}(\lambda_3)} \frac{(\alpha\lambda^2+\beta)e^{i\lambda x-\Omega(\lambda,t)}(g'_3\,\widetilde{)}(\Omega(\lambda,t),t)}{\Pi(\lambda)}d\lambda$$

$$+\int_{\mathcal{C}(\lambda_1)+\mathcal{C}(\lambda_2)+\mathcal{C}(\lambda_3)} \frac{e^{i\lambda x-\Omega(\lambda,t)}\widetilde{G}(\Omega(\lambda,t),\lambda,t)}{\Pi(\lambda)}d\lambda + \int_{-\infty}^{\infty} \frac{e^{i\lambda x-\Omega(\lambda,t)}\hat{f}(\Omega(\lambda,t),\lambda,t)}{\Pi(\lambda)}d\lambda ,$$

Next, we solve the algebraic equation

$$(8.7)\quad \frac{[\sigma(\lambda)]^2}{1+\gamma[\sigma(\lambda)]^2+\beta[\sigma(\lambda)]^4+\alpha[\sigma(\lambda)]^6} = \frac{\lambda^2}{1+\gamma\lambda^2+\beta\lambda^4+\alpha\lambda^6},$$

for $\sigma(\lambda)$. We may write (8.7) in the form

$$(8.8)\quad (\alpha\lambda^2)[\sigma(\lambda)]^6 + (\beta\lambda^2)[\sigma(\lambda)]^4 - (1+\beta\lambda^4+\alpha\lambda^6)[\sigma(\lambda)]^2 + \lambda^2,$$

which implies that the functions $\sigma(\lambda)$, which solve (8.7), are the square roots of the functions $\varsigma(\lambda)$ which satisfy the cubic equation

$$(8.9)\quad (\alpha\lambda^2)[\varsigma(\lambda)]^3 + (\beta\lambda^2)[\varsigma(\lambda)]^2 - (1+\beta\lambda^4+\alpha\lambda^6)[\varsigma(\lambda)]^2 + \lambda^2 = 0.$$

We procced, by fixing a $j \in \{1,2,3\}$ and considering the corresponding zero $\lambda_j$ of $\Pi(\lambda)$. We claim that there exists a biholomorphic function $\sigma_{j1} = \sigma_{j1}(\lambda)$, from an open neighborhood of $\lambda_j$ to an open neighborhood of $-\lambda_1$, such that $\sigma_{j1}(\lambda_j) = -\lambda_1$ and which satisfies (8.7), for $\lambda$ close to $\lambda_j$. Indeed, considering the polynomial

$$Q(\lambda,\sigma) = \alpha\lambda^2\sigma^6 + \beta\lambda^2\sigma^4 - (1+\beta\lambda^4+\alpha\lambda^6)\sigma^2 + \lambda^2 \in \mathbb{C}[\lambda,\sigma],$$

in the variables $\lambda$ and $\sigma$, we note that

$$(8.10)\quad Q(\lambda,\sigma)\big|_{(\lambda,\sigma)=(\lambda_j,-\lambda_1)} = 0 \quad\text{and}\quad \frac{\partial Q(\lambda,\sigma)}{\partial \sigma}\bigg|_{(\lambda,\sigma)=(\lambda_j,-\lambda_1)} \neq 0.$$

The second part of (8.10) follows from the fact that the polynomial $Q(\lambda_j,\sigma)$, considered as a polynomial of $\sigma$, has $\sigma = -\lambda_1$ as a simple zero. Therefore, the existence of $\sigma_{j1} = \sigma_{j1}(\lambda)$, with the required properties, follows from the implicit function theorem.

Similarly, there exist biholomorphic functions $\sigma_{j2} = \sigma_{j2}(\lambda)$ and $\sigma_{j3} = \sigma_{j3}(\lambda)$, from an open neighborhood of $\lambda_j$ to open neighborhoods of $-\lambda_2$ and $-\lambda_3$, respectively, such that $\sigma_{j2}(\lambda_j) = -\lambda_2$, $\sigma_{j3}(\lambda_j) = -\lambda_3$, and which satisfy (8.7), for $\lambda$ close to $\lambda_j$.

Using the functions $\sigma_{jk} = \sigma_{jk}(\lambda)$, $j,k \in \{1,2,3\}$, we obtain, in view of (8.5), the following equations:

$$(8.11)\quad \hat{u}(\sigma_{jk}(\lambda),t) = \hat{u}_0(\sigma_{jk}(\lambda)) - \frac{\alpha e^{-\Omega(\lambda,t)}(g'_5\,\widetilde{)}(\Omega(\lambda,t),t)}{\Pi(\sigma_{jk}(\lambda))}$$

$$-\frac{\alpha i\sigma_{jk}(\lambda)e^{-\Omega(\lambda,t)}(g'_4\,\widetilde{)}(\Omega(\lambda,t),t)}{\Pi(\sigma_{jk}(\lambda))} + \frac{(\alpha[\sigma_{jk}(\lambda)]^2+\beta)e^{-\Omega(\lambda,t)}(g'_3\,\widetilde{)}(\Omega(\lambda,t),t)}{\Pi(\sigma_{jk}(\lambda))}$$

$$+ \frac{e^{-\Omega(\lambda,t)}\widetilde{G}(\Omega(\lambda,t),\sigma_{jk}(\lambda),t)}{\Pi(\sigma_{jk}(\lambda))} + \frac{e^{-\Omega(\lambda,t)}\hat{f}(\Omega(\lambda,t),\sigma_{jk}(\lambda),t)}{\Pi(\sigma_{jk}(\lambda))}, \text{ for } \lambda \in \mathcal{C}(\lambda_j).$$

Fixing a $j \in \{1, 2, 3\}$, we consider the system of algebraic equations (8.11), with $k = 1, 2, 3$, and $\lambda \in \mathcal{C}(\lambda_j)$, in the unknown quantities

(8.12) $\qquad e^{-\Omega(\lambda,t)}(g'_5\widetilde{\ })(\Omega(\lambda,t),t), \; e^{-\Omega(\lambda,t)}(g'_4\widetilde{\ })(\Omega(\lambda,t),t), \; e^{-\Omega(\lambda,t)}(g'_3\widetilde{\ })(\Omega(\lambda,t),t), \; \lambda \in \mathcal{C}(\lambda_j).$

The determinant of the unknowns (8.12), in this system, is

$$\mathcal{D}(\lambda) = \det \begin{bmatrix} \dfrac{-\alpha}{\Pi(\sigma_{j,1}(\lambda))} & \dfrac{-\alpha i \sigma_{j,1}(\lambda)}{\Pi(\sigma_{j,1}(\lambda))} & \dfrac{\alpha[\sigma_{j,1}(\lambda)]^2 + \beta}{\Pi(\sigma_{j,1}(\lambda))} \\ \dfrac{-\alpha}{\Pi(\sigma_{j,2}(\lambda))} & \dfrac{-\alpha i \sigma_{j,2}(\lambda)}{\Pi(\sigma_{j,2}(\lambda))} & \dfrac{\alpha[\sigma_{j,2}(\lambda)]^2 + \beta}{\Pi(\sigma_{j,2}(\lambda))} \\ \dfrac{-\alpha}{\Pi(\sigma_{j,3}(\lambda))} & \dfrac{-\alpha i \sigma_{j,3}(\lambda)}{\Pi(\sigma_{j,3}(\lambda))} & \dfrac{\alpha[\sigma_{j,3}(\lambda)]^2 + \beta}{\Pi(\sigma_{j,3}(\lambda))} \end{bmatrix}.$$

A computation shows that the value of $\mathcal{D}(\lambda)$, at the point $\lambda = \lambda_j$, is equal to

(8.13) $\qquad \mathcal{D}(\lambda)\big|_{\lambda=\lambda_j} = \dfrac{\alpha^3 i \prod\limits_{1 \le k < l \le 3}[\sigma_{j,k}(\lambda_j) - \sigma_{j,l}(\lambda_j)]}{\Pi(\sigma_{j,1}(\lambda_j))\Pi(\sigma_{j,2}(\lambda_j))\Pi(\sigma_{j,3}(\lambda_j))} \ne 0,$

since, in view of (8.7), $\sigma_{j,1}(\lambda_j)$, $\sigma_{j,2}(\lambda_j)$, $\sigma_{j,3}(\lambda_j)$ are three of the six distinct zeros of the polynomial

$$\Pi(\sigma(\lambda)) = 1 + \gamma[\sigma(\lambda)]^2 + \beta[\sigma(\lambda)]^4 + \alpha[\sigma(\lambda)]^6.$$

It follows from (8.13) that $\mathcal{D}(\lambda) \ne 0$, for $\lambda$ close to $\lambda_1$, and, therefore, the quantities (8.12) can be computed from the system (8.11), with $k = 1, 2, 3$.

Another point, that we have to ensure and which is of utmost importance, is the vanishing of the integrals over the contours $\mathcal{C}(\lambda_j)$, which contain the terms $\hat{u}(\sigma_{jk}(\lambda),t)$, $j, k \in \{1, 2, 3\}$, and which arise when we substitute the quantities (8.12) in (8.6). But this is easy to prove, observing that the functions

$$\frac{1}{\Pi(\sigma_{j,k}(\lambda))\Pi(\sigma_{j,l}(\lambda))} \frac{1}{\mathcal{D}(\lambda)} \frac{1}{\Pi(\lambda)}, \; 1 \le k < l \le 3, \text{ are analytic for } \lambda \text{ close to } \lambda_j,$$

which – in turn – is a consequence of the analyticity of

$$\frac{\Pi(\sigma_{j,m}(\lambda))}{\Pi(\lambda)} = \frac{[\sigma_{j,m}(\lambda)]^2}{\lambda^2}, \; m \in \{1, 2, 3\}, \text{ for } \lambda \text{ close to } \lambda_j.$$

Thus, substituting in (8.6), the quantities (8.12), as these are computed by solving the systems (8.11), we express the solution $u(x,t)$, $x > 0$, $t > 0$, of problem (8.1), in terms of the data.

## 9. The operator $\partial_t - \alpha(t)\partial_{xxxxxt} + \beta(t)\partial_{xxxxx} - \gamma(t)$

**Problem 8** Solve

(9.1)
$$\begin{cases} \partial_t u = \alpha(t)\partial_{xxxxxt} u + \beta(t)\partial_{xxxxx} u - \gamma(t)u + f(x,t), \ (x,t) \in \mathbb{R}^+ \times \mathbb{R}^+, \\ u(x,0) = u_0(x), \ x \in \mathbb{R}^+, \\ u(0,t) = g_0(t), \ t \in \mathbb{R}^+, \\ u_x(0,t) = g_1(t), \ t \in \mathbb{R}^+, \\ u_{xx}(0,t) = g_2(t), \ t \in \mathbb{R}^+, \end{cases}$$

for $u(x,t)$.

**Derivation of the solution** For $\lambda \in \mathbb{C}$ with $\text{Im}\,\lambda \leq 0$, the differential equation in (9.1) leads to

$$[1+\alpha(t)\lambda^6]\frac{\partial}{\partial t}[\hat{u}(\lambda,t)] + [\beta(t)\lambda^6 + \gamma(t)]\hat{u}(\lambda,t)$$
$$= -\alpha(t)[g'_5(t) + i\lambda g'_4(t) + (i\lambda)^2 g'_3(t) + (i\lambda)^3 g'_2(t) + (i\lambda)^4 g'_1(t) + (i\lambda)^5 g'_0(t)]$$
$$- \beta(t)[g_5(t) + i\lambda g_4(t) + (i\lambda)^2 g_3(t) + (i\lambda)^3 g_2(t) + (i\lambda)^4 g_1(t) + (i\lambda)^5 g_0(t)] + \hat{f}(\lambda,t)$$

$\Rightarrow$

$$\frac{\partial}{\partial t}[\hat{u}(\lambda,t)] + \frac{\beta(t)\lambda^6 + \gamma(t)}{1+\alpha(t)\lambda^6}\hat{u}(\lambda,t)$$
$$= -\frac{\alpha(t)}{1+\alpha(t)\lambda^6}[g'_5(t) + i\lambda g'_4(t) - \lambda^2 g'_3(t) - i\lambda^3 g'_2(t) + \lambda^4 g'_1(t) + i\lambda^5 g'_0(t)]$$
$$- \frac{\beta(t)}{1+\alpha(t)\lambda^6}[g_5(t) + i\lambda g_4(t) - \lambda^2 g_3(t) - i\lambda^3 g_2(t) + \lambda^4 g_1(t) + i\lambda^5 g_0(t)] + \hat{f}(\lambda,t)$$

$\Rightarrow$

(9.2) $$\frac{\partial}{\partial t}[\hat{u}(\lambda,t)] + \frac{\beta(t)\lambda^6 + \gamma(t)}{1+\alpha(t)\lambda^6}\hat{u}(\lambda,t)$$
$$= -\frac{1}{1+\alpha(t)\lambda^6}[\alpha(t)g'_5(t) + \beta(t)g_5(t)] - \frac{i\lambda}{1+\alpha(t)\lambda^6}[\alpha(t)g'_4(t) + \beta(t)g_4(t)]$$
$$+ \frac{\lambda^2}{1+\alpha(t)\lambda^6}[\alpha(t)g'_3(t) + \beta(t)g_3(t)] + \frac{G(\lambda,t)}{1+\alpha(t)\lambda^6} + \frac{\hat{f}(\lambda,t)}{1+\alpha(t)\lambda^6},$$

where we have set

$$G(\lambda,t) = i\lambda^3[\alpha(t)g'_2(t) + \beta(t)g_2(t)] - \lambda^4[\alpha(t)g'_1(t) + \beta(t)g_1(t)] - i\lambda^5[\alpha(t)g'_0(t) + \beta(t)g_0(t)].$$

Furthermore, setting

$$G_5(t) = \alpha(t)g'_5(t) + \beta(t)g_5(t), \ G_4(t) = \alpha(t)g'_4(t) + \beta(t)g_4(t) \text{ and } G_3(t) = \alpha(t)g'_3(t) + \beta(t)g_3(t),$$

we write (9.2) as follows:

(9.3) $$\frac{\partial}{\partial t}[\hat{u}(\lambda,t)] + \frac{\beta(t)\lambda^6 + \gamma(t)}{1+\alpha(t)\lambda^6}\hat{u}(\lambda,t) = -\frac{G_5(t)}{1+\alpha(t)\lambda^6} - \frac{i\lambda G_4(t)}{1+\alpha(t)\lambda^6} + \frac{\lambda^2 G_3(t)}{1+\alpha(t)\lambda^6} + \frac{G(\lambda,t)}{1+\alpha(t)\lambda^6} + \frac{\hat{f}(\lambda,t)}{1+\alpha(t)\lambda^6}.$$

Defining

$$\omega(\lambda,t) = \frac{\beta(t)\lambda^6 + \gamma(t)}{1+\alpha(t)\lambda^6} \text{ and } \Omega(\lambda,t) = \int_{\tau=0}^{t}\omega(\lambda,\tau)d\tau = \int_{\tau=0}^{t}\frac{\beta(\tau)\lambda^6 + \gamma(\tau)}{1+\alpha(\tau)\lambda^6}d\tau,$$

we wriet (9.3) in thefollowing form

$$\hat{u}(\lambda,t) = e^{-\Omega(\lambda,t)}\hat{u}_0(\lambda) - e^{-\Omega(\lambda,t)}\int_{\tau=0}^{t}\frac{e^{\Omega(\lambda,\tau)}}{1+\alpha(\tau)\lambda^6}G_5(\tau)\,d\tau - e^{-\Omega(\lambda,t)}\int_{\tau=0}^{t}\frac{i\lambda e^{\Omega(\lambda,\tau)}}{1+\alpha(\tau)\lambda^6}G_4(\tau)\,d\tau$$

$$+ e^{-\Omega(\lambda,t)}\int_{\tau=0}^{t}\frac{\lambda^2 e^{\Omega(\lambda,\tau)}}{1+\alpha(\tau)\lambda^6}G_3(\tau)\,d\tau + e^{-\Omega(\lambda,t)}\int_{\tau=0}^{t}\frac{e^{\Omega(\lambda,\tau)}}{1+\alpha(\tau)\lambda^6}G(\lambda,\tau)\,d\tau + e^{-\Omega(\lambda,t)}\int_{\tau=0}^{t}\frac{e^{\Omega(\lambda,\tau)}}{1+\alpha(\tau)\lambda^6}\hat{f}(\lambda,\tau)\,d\tau,$$

i.e.,

(9.4) $\hat{u}(\lambda,t) = e^{-\Omega(\lambda,t)}\hat{u}_0(\lambda) - e^{-\Omega(\lambda,t)}\widetilde{G}_5(\lambda,t) - i\lambda e^{-\Omega(\lambda,t)}\widetilde{G}_4(\lambda,t) + \lambda^2 e^{-\Omega(\lambda,t)}\widetilde{G}_3(\lambda,t)$

$$+ e^{-\Omega(\lambda,t)}\widetilde{G}(\lambda,t) + e^{-\Omega(\lambda,t)}\widetilde{\hat{f}}(\lambda,t), \text{ for } \lambda \in \mathbb{C} \text{ with } \operatorname{Im}\lambda \leq 0,$$

where

$$\widetilde{G}_j(\lambda,t) := \int_{\tau=0}^{t}\frac{e^{\Omega(\lambda,\tau)}}{1+\alpha(\tau)\lambda^6}G_j(\tau)\,d\tau, \text{ for } j = 5,4,3,$$

$$\widetilde{G}(\lambda,t) := \int_{\tau=0}^{t}\frac{e^{\Omega(\lambda,\tau)}}{1+\alpha(\tau)\lambda^6}G(\lambda,\tau)\,d\tau \text{ and } \widetilde{\hat{f}}(\lambda,t) := \int_{\tau=0}^{t}\frac{e^{\Omega(\lambda,\tau)}}{1+\alpha(\tau)\lambda^6}\hat{f}(\lambda,\tau)\,d\tau.$$

Next, we fix a $t > 0$, and we find the zeros $\lambda_k = \lambda_k(t)$, $k = 1,2,3,4,5,6$, of the function $1+\alpha(t)\lambda^6$:

$$\lambda_1 = e^{i\pi/6}/\sqrt[6]{\alpha(t)},\ \lambda_2 = e^{i3\pi/6}/\sqrt[6]{\alpha(t)},\ \lambda_3 = e^{i5\pi/6}/\sqrt[6]{\alpha(t)},$$

$$\lambda_4 = -e^{i\pi/6}/\sqrt[6]{\alpha(t)},\ \lambda_5 = -e^{i3\pi/6}/\sqrt[6]{\alpha(t)},\ \lambda_6 = -e^{i5\pi/6}/\sqrt[6]{\alpha(t)}.$$

We also solve the equation

$$\omega(\lambda,t) = \omega(\sigma(\lambda),t) \Leftrightarrow \frac{\beta(t)[\sigma(\lambda)]^6 + \gamma(t)}{1+\alpha(t)[\sigma(\lambda)]^6} = \frac{\beta(t)\lambda^6 + \gamma(t)}{1+\alpha(t)\lambda^6},$$

and, for this, it suffices to have $[\sigma(\lambda)]^6 = \lambda^6$, i.e.,

$$\sigma_0(\lambda) = \lambda,\ \sigma_1(\lambda) = e^{i\pi/3}\lambda,\ \sigma_2(\lambda) = e^{i2\pi/3}\lambda,$$

$$\sigma_3(\lambda) = e^{i3\pi/3}\lambda = -\lambda,\ \sigma_4(\lambda) = e^{i4\pi/3}\lambda = -e^{i\pi/3}\lambda,\ \sigma_5(\lambda) = -e^{i2\pi/3}\lambda.$$

It is crucial for the construction that the above functions $\sigma_m(\lambda)$ satisfy the equations:

$$\Omega(\lambda,t) = \Omega(\sigma_m(\lambda),t) \text{ and } \widetilde{G}_j(\lambda,t) = \widetilde{G}_j(\sigma_m(\lambda),t),\ j = 5,4,3,\ m = 0,1,2,3,4,5.$$

Next, (9,4) gives

(9.6) $2\pi u(x,t) = \int_{-\infty}^{\infty}e^{i\lambda x}\hat{u}(\lambda,t)d\lambda$

$$= \int_{-\infty}^{\infty}e^{i\lambda x - \Omega(\lambda,t)}\hat{u}_0(\lambda)d\lambda - \int_{\partial A_1 + \partial A_2 + \partial A_3}e^{i\lambda x - \Omega(\lambda,t)}\widetilde{G}_5(\lambda,t)d\lambda - \int_{\partial A_1 + \partial A_2 + \partial A_3}i\lambda e^{i\lambda x - \Omega(\lambda,t)}\widetilde{G}_4(\lambda,t)d\lambda$$

$$+ \int_{\partial A_1 + \partial A_2 + \partial A_3}\lambda^2 e^{i\lambda x - \Omega(\lambda,t)}\widetilde{G}_3(\lambda,t)d\lambda + \int_{\partial A_1 + \partial A_2 + \partial A_3}e^{i\lambda x - \Omega(\lambda,t)}\widetilde{G}(\lambda,t)d\lambda + \int_{\partial A_1 + \partial A_2 + \partial A_3}e^{i\lambda x - \Omega(\lambda,t)}\widetilde{\hat{f}}(\lambda,t)d\lambda,\ x > 0,\ t > 0,$$

where $A_1,A_2,A_3$ are sufficiently small open neighborhoods of the half-lines

$$[c_0 e^{i\pi/6}, +\infty e^{i\pi/6}) := \{\eta e^{i\pi/6} : \eta \geq c_0\}, \ [c_0 e^{i3\pi/6}, +\infty e^{i3\pi/6}) := \{\eta e^{i3\pi/6} : \eta \geq c_0\} \text{ and}$$

$$[c_0 e^{i5\pi/6}, +\infty e^{i5\pi/6}) := \{\eta e^{i5\pi/6} : \eta \geq c_0\},$$

respectively, and whose boundaries $\partial \mathcal{A}_1, \partial \mathcal{A}_2, \partial \mathcal{A}_3$ are smooth and counterclockwise oriented. (We assume that $c_0 := \inf\{1/\sqrt[6]{\alpha(t)} : t \geq 0\} > 0$.)

Now, we have to compute the values of the integrals in (9.6), which contain the unknown quantities $\widetilde{G}_5(\lambda,t), \widetilde{G}_4(\lambda,t), \widetilde{G}_3(\lambda,t)$, in terms of the data. This will be done by solving three systems of algebraic equations – each system will be solved for $\lambda \in \partial \mathcal{A}_1, \lambda \in \partial \mathcal{A}_2, \lambda \in \partial \mathcal{A}_3$, separately.

Thus, for $\lambda \in \partial \mathcal{A}_1$, we have, in view of (9.4),

(9.7) $\hat{u}(-\lambda,t) = e^{-\Omega(\lambda,t)}\hat{u}_0(-\lambda) - e^{-\Omega(\lambda,t)}\widetilde{G}_5(\lambda,t) + i\lambda e^{-\Omega(\lambda,t)}\widetilde{G}_4(\lambda,t) + \lambda^2 e^{-\Omega(\lambda,t)}\widetilde{G}_3(\lambda,t)$
$$+ e^{-\Omega(\lambda,t)}\widetilde{G}(-\lambda,t) + e^{-\Omega(\lambda,t)}\widetilde{\tilde{f}}(-\lambda,t),$$

(9.8) $\hat{u}(\sigma_4(\lambda),t) = e^{-\Omega(\lambda,t)}\hat{u}_0(\sigma_4(\lambda)) - e^{-\Omega(\lambda,t)}\widetilde{G}_5(\lambda,t) - i\sigma_4(\lambda)e^{-\Omega(\lambda,t)}\widetilde{G}_4(\lambda,t) + [\sigma_4(\lambda)]^2 e^{-\Omega(\lambda,t)}\widetilde{G}_3(\lambda,t)$
$$+ e^{-\Omega(\lambda,t)}\widetilde{G}(\sigma_4(\lambda),t) + e^{-\Omega(\lambda,t)}\widetilde{\tilde{f}}(\sigma_4(\lambda),t),$$

(9.9) $\hat{u}(\sigma_5(\lambda),t) = e^{-\Omega(\lambda,t)}\hat{u}_0(\sigma_5(\lambda)) - e^{-\Omega(\lambda,t)}\widetilde{G}_5(\lambda,t) - i\sigma_5(\lambda)e^{-\Omega(\lambda,t)}\widetilde{G}_4(\lambda,t) + [\sigma_5(\lambda)]^2 e^{-\Omega(\lambda,t)}\widetilde{G}_3(\lambda,t)$
$$+ e^{-\Omega(\lambda,t)}\widetilde{G}(\sigma_5(\lambda),t) + e^{-\Omega(\lambda,t)}\widetilde{\tilde{f}}(\sigma_5(\lambda),t).$$

Next, we solve the above system for the unknown quantities

(9.10) $\qquad e^{-\Omega(\lambda,t)}\widetilde{G}_5(\lambda,t), \ e^{-\Omega(\lambda,t)}\widetilde{G}_4(\lambda,t), \ e^{-\Omega(\lambda,t)}\widetilde{G}_3(\lambda,t)$, for $\lambda \in \partial \mathcal{A}_1$.

It is crucial to ensure that the determinant of the unknowns (9.10) is different from zero. Indeed,

$$\det \begin{bmatrix} -1 & i\lambda & \lambda^2 \\ -1 & -i\sigma_4(\lambda) & [\sigma_4(\lambda)]^2 \\ -1 & -i\sigma_5(\lambda) & [\sigma_5(\lambda)]^2 \end{bmatrix} \neq 0,$$

since the numbers $1, e^{-i\pi/3}, e^{-i2\pi/3}$ are three of the six distinct 6th roots of unity.

Similarly, we write systems analogous to (9.7), (9.8) and (9.9), for the unknowns quantities

(9.11) $\qquad e^{-\Omega(\lambda,t)}\widetilde{G}_5(\lambda,t), \ e^{-\Omega(\lambda,t)}\widetilde{G}_4(\lambda,t), \ e^{-\Omega(\lambda,t)}\widetilde{G}_3(\lambda,t)$, for $\lambda \in \partial \mathcal{A}_2$,

and

(9.10) $\qquad e^{-\Omega(\lambda,t)}\widetilde{G}_5(\lambda,t), \ e^{-\Omega(\lambda,t)}\widetilde{G}_4(\lambda,t), \ e^{-\Omega(\lambda,t)}\widetilde{G}_3(\lambda,t)$, for $\lambda \in \partial \mathcal{A}_3$,

and solve them.

Finally, substituting the values of (9.10), (9.11) and (9.12), in (9.6), we obtain the integral representation of the solution $u(x,t)$, in terms of the data. As usual, it is crucial to ensure that the integrals, which arise, in (9.6), by carrying out this substitution and which contain the unknown quantities $\hat{u}(-\lambda,t), \hat{u}(\sigma_4(\lambda),t), \hat{u}(\sigma_5(\lambda),t)$, vanish.

## 10. The operator $\partial_t + (-1)^\nu \alpha(t)\partial_x^{2\nu}\partial_t + (-1)^\nu \beta(t)\partial_x^{2\nu} + \gamma(t)$

**Problem 9** *Solve*

(10.1)
$$\begin{cases} \partial_t u = -(-1)^\nu \alpha(t)\partial_x^{2\nu}\partial_t u - (-1)^\nu \beta(t)\partial_x^{2\nu} u - \gamma(t)u + f(x,t), \ (x,t)\in\mathbb{R}^+\times\mathbb{R}^+, \\ u(x,0) = u_0(x), \ x\in\mathbb{R}^+, \\ \partial_x^{j-1} u(0,t) = g_{j-1}(t), \ j=1,2,...,\nu, \ t\in\mathbb{R}^+, \end{cases}$$

*for* $u(x,t)$.

***Outline of construction of the solution*** In light of the formula

$$\int_0^\infty \frac{\partial^{2\nu} u(x,t)}{\partial x^{2\nu}} e^{-i\lambda x} dx = -\sum_{k=1}^{2\nu}(i\lambda)^{k-1} g_{2\nu-k}(t) + (i\lambda)^{2\nu}\hat{u}(\lambda,t),$$

where

$$g_{2\nu-k}(t) = \left.\frac{\partial^{2\nu-k} u(x,t)}{\partial x^{2\nu}}\right|_{x=0},$$

and the previously analysed cases, i.e., when $\nu=2$ (in section 2), $\nu=4$ (in section 4), and $\nu=6$ (in section 9), the differential equation in (10.1) leads to

(10.2) $\dfrac{\partial}{\partial t}[\hat{u}(\lambda,t)] + \dfrac{\beta(t)\lambda^{2\nu}+\gamma(t)}{1+\alpha(t)\lambda^{2\nu}}\hat{u}(\lambda,t)$

$$= -\frac{\alpha(t)}{1+\alpha(t)\lambda^{2\nu}}\sum_{k=1}^{2\nu}(i\lambda)^{k-1} g'_{2\nu-k}(t) - \frac{\beta(t)}{1+\alpha(t)\lambda^{2\nu}}\sum_{k=1}^{2\nu}(i\lambda)^{k-1} g_{2\nu-k}(t) + \frac{\hat{f}(\lambda,t)}{1+\alpha(t)\lambda^6}$$

$$= -\frac{1}{1+\alpha(t)\lambda^{2\nu}}\sum_{k=1}^{2\nu}(i\lambda)^{k-1}[\alpha(t)g'_{2\nu-k}(t) + \beta(t)g_{2\nu-k}(t)] + \frac{\hat{f}(\lambda,t)}{1+\alpha(t)\lambda^6}.$$

Defining

$$\Omega(\tau) = \int_{\tau=0}^t \frac{\beta(\tau)\lambda^{2\nu}+\gamma(\tau)}{1+\alpha(\tau)\lambda^{2\nu}} d\tau,$$

$$G_{2\nu-k}(t) = \alpha(t)g'_{2\nu-k}(t) + \beta(t)g_{2\nu-k}(t), \text{ for } k=1,2,...,\nu,$$

and

$$G(\lambda,t) = -\sum_{k=\nu+1}^{2\nu}(i\lambda)^{k-1}[\alpha(t)g'_{2\nu-k}(t) + \beta(t)g_{2\nu-k}(t)],$$

we obtain, in view of (10.2), for $\lambda\in\mathbb{C}$ with $\operatorname{Im}\lambda \leq 0$,

(10.3) $\hat{u}(\lambda,t) = e^{-\Omega(\lambda,t)}\hat{u}_0(\lambda) - \sum_{k=1}^\nu (i\lambda)^{k-1}e^{-\Omega(\lambda,t)}\widetilde{G}_{2\nu-k}(\lambda,t) + e^{-\Omega(\lambda,t)}\widetilde{G}(\lambda,t) + e^{-\Omega(\lambda,t)}\widetilde{\hat{f}}(\lambda,t),$

where we have set

$$\widetilde{G}_{2\nu-k}(\lambda,t) := \int_{\tau=0}^t \frac{e^{\Omega(\lambda,\tau)}}{1+\alpha(\tau)\lambda^{2\nu}} G_{2\nu-k}(\tau)\, d\tau, \text{ for } k=1,2,...,\nu,$$

$$\widetilde{G}(\lambda,t) := \int_{\tau=0}^t \frac{e^{\Omega(\lambda,\tau)}}{1+\alpha(\tau)\lambda^{2\nu}} G(\lambda,\tau)\, d\tau \text{ and } \widetilde{\hat{f}}(\lambda,t) := \int_{\tau=0}^t \frac{e^{\Omega(\lambda,\tau)}}{1+\alpha(\tau)\lambda^{2\nu}} \hat{f}(\lambda,\tau)\, d\tau.$$

Next, we fix a $t > 0$ and we find the zeros $\lambda_j = \lambda_j(t)$, $1 \leq j \leq \nu$, of the function $1 + \alpha(t)\lambda^{2\nu}$, which lie in the upper half plane, in which case the numbers $-\lambda_j = -\lambda_j(t)$, $1 \leq j \leq \nu$, are the zeros of $1 + \alpha(t)\lambda^{2\nu}$ in the lower half plane. Then, (10.3) leads to the following representation formula:

$$(10.4) \quad 2\pi u(x,t) = \sum_{j=1}^{\nu} \int_{\mathcal{A}_j} e^{i\lambda x - \Omega(\lambda,t)} \hat{u}_0(\lambda) d\lambda - \sum_{\kappa=1}^{\nu} \sum_{j=1}^{\nu} \int_{\mathcal{A}_j} (i\lambda)^{k-1} e^{i\lambda x - \Omega(\lambda,t)} \widetilde{G}_{2\nu-k}(\lambda,t) d\lambda$$

$$+ \sum_{k=1}^{\nu} \sum_{j=1}^{\nu} \int_{\mathcal{A}_j} e^{i\lambda x - \Omega(\lambda,t)} \widetilde{G}(\lambda,t) d\lambda + \sum_{k=1}^{\nu} \sum_{j=1}^{\nu} \int_{\mathcal{A}_j} e^{i\lambda x - \Omega(\lambda,t)} \widetilde{\hat{f}}(\lambda,t) d\lambda, \quad x > 0, \, t > 0,$$

where $\mathcal{A}_1, \mathcal{A}_2, \ldots, \mathcal{A}_\nu$ are sufficiently and appropriately small open neighborhoods of the half-lines

$$[c_0 e^{i\pi[1+2(j-1)]/2\nu}, +\infty e^{i\pi[1+2(j-1)]/2\nu}) := \{\eta e^{i\pi[1+2(j-1)]/2\nu} : \eta \geq c_0\}, \quad 1 \leq j \leq \nu,$$

respectively, and whose boundaries $\partial \mathcal{A}_1, \partial \mathcal{A}_2, \ldots, \partial \mathcal{A}_\nu$ are smooth and counterclockwise oriented. (We assume that $c_0 := \inf\{1/\sqrt[2\nu]{\alpha(t)} : t \geq 0\} > 0$.)

Finally, in order to compute the integrals in (10.4), which contain the unknown quantities

$$e^{-\Omega(\lambda,t)} \widetilde{G}_{2\nu-k}(\lambda,t), \text{ for } \lambda \in \mathcal{A}_j \text{ with } 1 \leq j \leq \nu \text{ and for each } 1 \leq k \leq \nu,$$

we set up $\nu$ systems of $\nu$ algebraic equations in $\nu$ unknowns each, with each system corresponding to each contour $\mathcal{A}_j$, using the equation (10.3) and the appropriate symmetry relations, thereby obtaining an effective UTM solution formula for the "canonical" problem 9.